%% file: NonparChoice.tex
\documentclass[opre,nonblindrev]{informs3mod}

\OneAndAHalfSpacedXII %

\input{Def}

\usepackage{algorithm,algorithmic}

\usepackage[colorlinks=true,breaklinks=true,bookmarks=true,urlcolor=blue,citecolor=blue,linkcolor=blue,bookmarksopen=false,draft=false]{hyperref}

\usepackage{caption,subcaption}
\usepackage{bm,enumerate,verbatim}
\usepackage{multirow}
\usepackage{tikz}
\usepackage{pgfplots}
\pgfplotsset{compat=1.3}

\theoremstyle{example}

\DeclareMathOperator{\train}{train}
\DeclareMathOperator{\test}{test}
\DeclareMathOperator{\MAE}{MAE}

\usepackage{natbib}
 \bibpunct[, ]{(}{)}{,}{a}{}{,}%
 \def\bibfont{\small}%
 \def\BIBand{and}%
 
\usepackage{multibib}
\newcites{EC}{References for Appendix}

\TheoremsNumberedThrough     %

\EquationsNumberedThrough    %

\MANUSCRIPTNO{} %

\begin{document}

\RUNAUTHOR{Ho-Nguyen and K{\i}l{\i}n\c{c}-Karzan}

\RUNTITLE{Dynamic Data-Driven Estimation of Non-Parametric Choice Models}

\TITLE{Dynamic Data-Driven Estimation of Non-Parametric Choice Models}

\ARTICLEAUTHORS{%
\AUTHOR{Nam Ho-Nguyen}
\AFF{Discipline of Business Analytics, The University of Sydney}
\AUTHOR{Fatma K{\i}l{\i}n\c{c}-Karzan}
\AFF{Tepper School of Business, Carnegie Mellon University} %
} %

\ABSTRACT{%
We study non-parametric estimation of choice models, which were introduced to alleviate unreasonable assumptions in traditional parametric models, and are prevalent in several application areas. Existing literature focuses only on the static observational setting where all of the observations are given upfront, they are not equipped with explicit convergence rate guarantees, and consequently they cannot provide an a priori analysis for the model accuracy vs sparsity trade-off on the actual estimated model returned by their algorithms. As opposed to this, we focus on estimating a non-parametric choice model from  observational data in a \emph{dynamic} setting, where observations are obtained over time. We show that choice model estimation can be cast as a convex-concave saddle-point (SP) joint estimation and optimization (JEO) problem, and we provide a primal-dual framework for deriving algorithms to solve this based on online convex optimization. By tailoring our framework carefully to the choice model estimation problem, we obtain tractable algorithms with provable convergence guarantees and explicit bounds on the sparsity of the estimated model. Our numerical experiments confirm the effectiveness of the algorithms derived from our framework.
}%


\maketitle

\input{Intro}

\input{DynamicLearning-choice}

\input{DynamicLearning}

\input{Algorithms}

\input{ExperimentsMain}

\input{Conclusion}

\ACKNOWLEDGMENT{This research is supported in part by NSF grant CMMI 1454548, NSF Award 1740707, and Award N660011824020 from the DARPA Lagrange Program.}

\input{NonparChoice.bbl}
\ECSwitch
\ECHead{Appendices to \emph{Dynamic Data-Driven Estimation of Non-Parametric Choice Models}}

\begin{APPENDICES}
\input{Proofs}

\input{ExistingApproaches}

\input{ErrorRates}

\input{OCOframework}

\input{FOMComparisons}

\input{FWcomparison}

\input{Experiments}

\input{EC.bbl}
\end{APPENDICES}

\end{document}

%% file: Def.tex
\definecolor{gold}{rgb}{0.85,0.65,0}

\renewcommand{\theequation}{\mbox{\arabic{equation}}}
\newcommand{\ese}{\end{eqnarray*}}
\newcommand{\bse}{\begin{eqnarray*}}
\def\beq{\begin{equation}}
\def\eeq{\end{equation}}

\def\fnote#1{\footnote}
\newcommand{\epr}{\hfill\hbox{\hskip 4pt \vrule width 5pt height 6pt depth 1.5pt}\vspace{0.0cm}\par}

\def\ra{\rangle}
\def\la{\langle}

\newcommand{\grad}{\ensuremath{\nabla}}

\def\I{{\mathbb{I}}}

\def\N{{\mathbb{N}}}

\def\R{{\mathbb{R}}}

\def\cA{\mathcal{A}}

\def\cR{{\cal R}}
\def\cS{\mathcal{S}}

\newcommand{\bbE}{\mathbb{E}}

\newcommand{\bbN}{\mathbb{N}}

\newcommand{\bbP}{\mathbb{P}}

\newcommand{\bbR}{\mathbb{R}}

\newcommand{\CA}{\mathcal{A}}

\newcommand{\CD}{\mathcal{D}}

\newcommand{\CP}{\mathcal{P}}

\newcommand{\CR}{\mathcal{R}}

\DeclareMathOperator{\SV}{SV}
\DeclareMathOperator{\sad}{sad}

\def\epsilonsad{\epsilon_{\sad}}

\DeclareMathOperator{\KL}{KL}

\def\dim{\mathop{{\rm dim}\,}}

\DeclareMathOperator{\Conv}{Conv}
\DeclareMathOperator{\conv}{conv}

\def\log{\mathop{{\rm log}}}

%% file: Intro.tex
\section{Introduction}\label{sec:intro}

A choice model is an effective tool to summarize and understand the preferences of a population over a set of items.
Such models give choice probabilities, that is, the probability that an agent will choose a particular item from a given subset.
They are prevalent in several application areas such as revenue management, web page ranking, betting theory, social choice, marketing, and economics (see \citet{DworkKumarNaorSivakumar2001,TallurivanRyzin2006book,JagabathulaShah2008,FariasJagShah2009,DesirGoyalJagSegev2016} and references therein). A good choice model aims to capture complex substitution behaviors of agents in order to accurately describe preferences from limited observations.

Choice model estimation has received quite a bit of interest. Traditional choice models often specify a \emph{parametric} structure for the choice probabilities (examples include the multinomial logit (MNL),
nested logit, and mixed MNL models); see \citet{TallurivanRyzin2006book}
and references therein.
Imposing a parametric structure makes estimation of the necessary parameters a simpler task, but is often at the expense of 
overly facile assumptions on substitution behaviors (such as independence of irrelevant alternatives in MNL models) and consequently preventing us from
accurately capturing preferences.
Therefore, the \emph{non-parametric} approach of directly estimating a probability distribution over rankings has drawn growing interest in academia and in practice \citep{Rusmevichientong2006,FariasJS2013,FariasJS2017},  
and through case studies, it is shown to lead to substantial improvement in prediction accuracy (see \citet{FariasJS2013,JagabathulaRusmevichientong2016}). 
In this paper, we focus on the problem of estimating a non-parametric choice model in a dynamic observation setting, and develop convex optimization-based approaches equipped with convergence and sparsity guarantees.

\textbf{Related Literature.} Earliest studies on non-parametric choice models appear in the economics and psychology literatures, e.g., \citet{BlockMarschak1960}. \citet{MahajanVanRyzin2001} showed that non-parametric models capture a number of parametric models as special cases.
The recent literature focuses on the \emph{static} estimation of non-parametric choice models, %
where the aim is to find a model that either matches the observed empirical probabilities with the model-based choice probabilities exactly (see \citet{FariasJS2013} which uses a dual-based constraint sampling approach), or uses column generation to minimize a \emph{distance measure} between the two based on Kullback-Leibler (KL) divergence \citep{vanRyzinVulcano2015} or $\ell_1$-norm \citep{BertsimasMisic2015}.

As opposed to the specific distance measures used in the prior literature, the recent work of \citet{JagabathulaRusmevichientong2016} %
focuses on general distance measures. While they suggest using the Frank-Wolfe (F-W) algorithm to estimate a non-parametric choice model, the bulk of their work focuses on a particular combinatorial subproblem that arises in all of the methods for non-parametric choice model estimation. 
Their main contribution is the characterization of sufficient conditions (in terms of the subset structure of the items) under which this subproblem becomes polynomial-time solvable.

In the static setting, these four approaches, \citet{FariasJS2013, vanRyzinVulcano2015, BertsimasMisic2015}, and  \citet{JagabathulaRusmevichientong2016}, are closely related to our work, thus we further discuss and compare them in Appendix~\ref{sec:learning-choice-model-existing}. We now discuss a number of important considerations in the estimation of non-parametric choice models: convergence of the algorithm, sparsity,
and dynamic data.

\paragraph{Convergence.} \citet{vanRyzinVulcano2015,BertsimasMisic2015} both suggest column generation procedures (for different objectives) that terminate in finite time. However, convergence \emph{rates} (i.e., how many iterations needed to obtain an $\epsilon$-accurate model) are not given for their column generation procedure. \citet{JagabathulaRusmevichientong2016} proposed using the F-W algorithm, for which convergence rates exist. Nevertheless, a direct application of F-W method to classical distance measures such as norms used in choice modeling is not possible due to the non-smoothness of these measures. \citet{FariasJS2013} does not discuss convergence of their method.

\paragraph{Sparsity.} In a full non-parametric model, there are a factorial number of probabilities to estimate, and thus even specifying a full non-parametric model is intractable even for moderate-sized problems, let alone computing choice probabilities or estimating one from data. This necessarily places an importance on characterizing the model simplicity (sparsity) vs model accuracy trade-off. \citet[Theorem 4]{FariasJS2013} provide an existence guarantee of some sparse model that fits the data. But, they do not provide sparsity guarantees on the \emph{actual} estimated model obtained from their proposed algorithm. The other three approaches aim for model sparsity by increasing the support of the distribution by at most one at each iteration. However, since none of them are equipped with explicit convergence rate guarantees,
they cannot provide an a priori analysis for the model accuracy vs sparsity trade-off on the actual estimated model returned by their algorithms.

\paragraph{Dynamic data.} Existing techniques for estimation of non-parametric choice models are not designed to efficiently work with the dynamic data, i.e., exploit the possibility of continuously updating the empirical choice probabilities as more observations are collected.
This setup is very realistic with today's data collection capabilities.
A na\"ive way to work with dynamic data is to simply re-solve the estimation problem each time we update the choice probabilities. However, in the case of estimation of non-parametric choice models, the static estimation problem is already very expensive to solve, thus such a na\"ive approach significantly compounds the existing computational challenges.

The convex optimization literature has recently considered the problem of minimizing a function given only an approximate sequence converging to true parameters; referred to as \emph{misspecified optimization} \citep{AhmadiShanbhag2014,JiangShanbhag2014} or \emph{joint estimation-optimization} (JEO, \citet{HoNguyenKK2016}). These methods provide a way to solve this problem without re-optimizing each time a new set of parameters is given, but instead performing only a single projection-type operation at each time step. The dynamic data setting for choice model estimation exactly fits this framework, where we think of the approximate sequence as the empirical choice probabilities computed over a growing set of observational data. Under mild statistical assumptions, the empirical probabilities will converge almost surely. However, the existing methods for JEO are ill-equipped to deal with the challenging non-parametric choice model estimation problem, since the domain of interest in general only admits a high-dimensional representation over which projection-type operations are difficult to perform.

\textbf{Contributions and Outline.} In this paper, we simultaneously address the points raised above on convergence, sparsity and dynamic data for non-parametric choice model estimation. We first describe in Section \ref{sec:learning-choice-model} that the dynamic non-parametric choice model estimation problem, and show that it can be solved as a dynamic JEO saddle point (SP) problem. Inspired by the distance measures used in the prior literature for static non-parametric choice estimation, we examine a large class of distance measures and show that they admit favorable SP structure. Importantly, this class includes all norms, but also covers well-known smoothed versions of these norms as well.

We then present a primal-dual framework for solving the JEO SP problem in Section \ref{sec:JEO}. Specifically, we provide a decomposition of the SP gap into a regret term and two error terms, and describe how the error terms behave as our data approximations become more accurate. In Section \ref{sec:JEO-algs}, we describe regret minimization, which can be used to bound the regret term in the saddle point gap decomposition.

We suggest in Section~\ref{sec:regret-algs} three possible regret minimizing algorithms, and give their regret bounds explicitly. Combined with the error terms, the regret bounds provide rigorous convergence guarantees for our framework; see Remark \ref{rem:rates}. Our convergence rates also imply explicit guarantees on the sparsity of our estimated choice model (see Remarks \ref{rem:sparsity-guarantee} and \ref{rem:rates}). Consequently, our results highlight a natural trade-off between desired estimation accuracy and model sparsity.

The high-dimensional structure of the domain of the choice model estimation problem necessitated us to extend the existing JEO literature to the SP setting; see Remark \ref{rem:linear-opt-X}. By viewing it as a SP problem and employing our primal-dual framework, we avoid the problem of having to solve difficult non-linear problems on the high-dimensional domain, and instead solve a relatively compact integer linear program (see Section \ref{sec:comb-subproblem}) at each iteration, while still enjoying efficient convergence guarantees. Furthermore, we believe that this JEO SP framework can be of interest beyond choice model estimation.

Our JEO SP framework has ample flexibility that it leads to a variety of solution algorithms simply by utilizing different regret minimization algorithms.
In particular, a slight variant of the Frank-Wolfe algorithm (also known as the Conditional Gradient algorithm) can be derived from our framework.  Our JEO SP framework in addition allows us to provide a new regularity condition to ensure convergence and avoid accumulation of errors (from using inexact data) which is different to the approximate gradient assumption
seen in the usual F-W literature; see Remarks \ref{rem:comparison-FW} and \ref{rem:interpretations}.

It is well-known that the F-W method in general \emph{does not} converge on non-smooth objectives; see \citet[Example 1]{Nesterov2018}, and we show in Appendix~\ref{sec:FWcomparison} that the classical distance measures used in the non-parametric choice estimation setup, such as norms and KL divergence are non-smooth. Our framework allows us to also derive algorithms for the non-smooth setting by selecting appropriate regret minimization algorithms, something which the usual F-W algorithm is unable to do. 
In Section \ref{sec:computationSummary}, we carry out a numerical study and compare several algorithms derived from our framework with some existing methods.  We find that the methods from our framework working directly with the non-smooth distance measures outperform the ones that require smooth approximations, including  the existing ones from the literature. We give brief concluding remarks in Section~\ref{sec:conclusion}.

We present all of the proofs in Appendix~\ref{sec:proofs}. Existing approaches to  non-parametric choice estimation are discussed in detail in Appendix~\ref{sec:learning-choice-model-existing}. Appendix \ref{sec:error-rates} provides rates for certain error terms in our analysis. Appendix \ref{sec:oco-intro} introduces online convex optimization, a key concept in our framework.  Appendix \ref{sec:interpretations} relates some algorithms derived from our framework to existing ones, such as Mirror Descent and the F-W algorithm. Appendix~\ref{sec:FWcomparison}, discusses different approaches for smoothing norm-based distance measures and compares our methods with existing F-W based approaches. Full details of our experimental setup from Section \ref{sec:computationSummary} and supplementary numerical results are given in Appendix~\ref{sec:computationApp}.

\textbf{Notation.} For a positive integer $n\in\N$,  we let $[n]:=\{1,\ldots,n\}$, define $\Delta_n:=\{x\in\R^n_+:~\sum_{i\in[n]} x_i=1\}$ to be the standard simplex, and $S_n$ to be the collection of rankings/permutations of the set $[n]$.
We denote the %
identity matrix in $\bbR^{n\times n}$ by $I_n$. 
We refer to a collection of objects $b_j$, $j \in J$ by the notation $\{b_j\}_{j \in J}$.
Given vectors $x$ and $y$, $\la x,y\ra$ corresponds to the usual inner product of $x$ and $y$. Given a norm $\|\cdot\|$ on a Euclidean space $\bbE$, we denote its dual norm by $\|x\|_* := \min_{y} \{ \la x, y\ra:~\|x\| \leq1 \}$.
For $q\in[1,\infty]$, $\|x\|_q$ denotes the usual $\ell_q$ norm of $x$. 
For a convex function $f$, we abuse notation slightly by denoting $\grad f(x)$ for both the gradient of function $f$ at $x$ if $f$ is differentiable and a subgradient of $f$ at $x$, even if $f$ is not differentiable. If $\phi$ is of the form $\phi(x,y)$, then $\grad_x \phi(x,y)$ denotes the (sub)gradient of $\phi$ at $x$ while keeping the other variables fixed at $y$.

%% file: DynamicLearning-choice.tex
\section{Dynamic Estimation of a Non-Parametric Choice Model}\label{sec:learning-choice-model}

\subsection{Model and Data}\label{sec:model-data}

In choice modeling, given a set of $n$ items, $[n] = \{1,\ldots,n\}$, we would like to understand the preferences of a population by estimating a model from choice observations. We follow a non-parametric approach to choice modeling introduced by \citet{FariasJS2013}. 
A non-parametric choice model is described by a probability distribution $\lambda \in \Delta_{n!}$ over all rankings $S_n$ of the items $[n]$. Given a ranking $\sigma \in S_n$, we think of $\lambda(\sigma)$ as the probability that a member of the population will rank the items according to $\sigma$. Also, when a particular member with ranking $\sigma$ is presented a subset of items $A \subseteq [n]$, they will choose the highest $\sigma$-ranked item $i(\sigma,A) = \argmin_{i \in A} \sigma(i)$. Thus, the probability of a random member of the population choosing an item $i \in A$ when presented with a subset $A$ is
\vspace{-7.5pt}
\[ \bbP_{\lambda}[i \mid A] := \sum_{\sigma \in S_n(i,A)} \lambda(\sigma), \quad S_n(i,A) := \{\sigma \in S_n : \text{$i$ is the highest $\sigma$-ranked item in $A$}\}. \]

Our choice observation set can be described as a collection of $K$ pairs $\left\{i^k,A^k\right\}_{k=1}^K$, where $i^k \in  A^k$ is the item chosen when the subset of items $A_k \subseteq [n]$ was presented. There are a finite number of possible subsets amongst the observations, we denote these by $A_j$, $j \in [m] = \{1,\ldots,m\}$. We also denote $N := \sum_{j=1}^{m} |A_j|$. In practice, the collection of possible subsets $\{A_j\}_{j \in [m]}$ can be controlled. Indeed, structural properties of these can have an impact on a combinatorial subproblem that appears in all non-parametric choice model estimation methods (see Section \ref{sec:comb-subproblem}). However, in this paper, we take this collection as given; see \citet{JagabathulaRusmevichientong2016} for a study on how the structure of $\{A_j\}_{j \in [m]}$ impacts the combinatorial subproblem. Based on this observation set, we define
\vspace{-7.5pt}
\begin{equation}
q_{ij} := \frac{1}{K} \sum_{k=1}^{K} \I(i^k = i, A^k = A_j), \quad q_{j} := \frac{1}{K} \sum_{k=1}^{K} \I(A^k = A_j) \quad\text{and}\quad p_{ij} := \frac{q_{ij}}{q_j},\label{eqn:pq_{ij}}
\end{equation}
where $\I$ denotes the indicator function, i.e., $\I(\cS) = 1$ if statement $\cS$ holds, and $\I(\cS) = 0$ otherwise. 
In words, $q_{ij}$ is the proportion of observations where assortment $A_j$ was displayed and item $i$ was chosen, $q_j$ is the proportion of observations where assortment $A_j$ was displayed, and $p_{ij}$ is the proportion of observations where item $i$ was chosen amongst those where assortment $A_j$ was displayed. Indeed, $p_{ij}$ are the empirical choice probabilities which we will use to tune the probability distribution $\lambda$. We denote the collection of these empirical choice probabilities as $p = \{p_{ij}\}_{i \in A_j, j \in [m]} \in \bbR^N$.

Our goal is to tune $\lambda$ so that the choice probabilities $\bbP_{\lambda}[i \mid A_j]$ are close to the empirical probabilities $p_{ij}$. We now define notation to succinctly describe $\bbP_{\lambda}[i \mid A_j]$. For a given pair $i \in A_j$ and ranking $\sigma \in S_n$, we let $a_{ij}(\sigma) = 1$ if $i$ is the highest $\sigma$-ranked item in $A_j$ (i.e., $\sigma \in S_n(i,A_j)$), and $a_{ij}(\sigma) = 0$ otherwise. We define vectors $a_{ij} = \{a_{ij}(\sigma)\}_{\sigma \in S_n} \in \{0,1\}^{n!}$ and $a(\sigma) = \{a_{ij}(\sigma)\}_{i \in A_j, j \in [m]} \in \{0,1\}^N$. Then 
$%
\bbP_{\lambda}[i \mid A_j] = \sum_{\sigma \in S_n} a_{ij}(\sigma) \lambda(\sigma) = \la a_{ij}, \lambda \ra. $ %

The polytope $X$ of all possible choice probabilities on observed pairs $i \in A_j$ consistent with some distribution $\lambda$ is
\vspace{-10pt}
\begin{equation}\label{eqn:X-choice-domain}
X := \conv\left( \left\{a(\sigma) :~ \sigma \in S_n \right\} \right) \subseteq \bbR^N.
\end{equation}
Our goal can now be stated as finding a point $x \in X$ that is close to $p$, and our choice model will be the weights $\lambda$ such that $x = \sum_{\sigma \in S_n} \lambda(\sigma) a(\sigma)$. Formally, we solve
\vspace{-7.5pt}
\begin{equation}\label{eqn:general-distance-learn-choice}
\min_x \left\{ D(x,p):~ x \in X \right\}
\end{equation}
for some appropriately defined distance measure $D(x,p)$, e.g., $D(x,p)=\|x-p\|$.

Finally, we describe the dynamic data setting, where we obtain additional observations over time. We denote a point in time by $t \in \bbN$, and the number of observations collected by time $t$ as $K_t \in \bbN$, which is non-decreasing in $t$. The set of observations at time $t$ is $\{i^k,A^k\}_{k \in [K_t]}$. For simplicity, we assume the collection of observed subsets $\{A_j\}_{j \in [m]}$ remains the same over time. We can compute empirical probabilities $p_{ij}^t$ at time $t$ by \eqref{eqn:pq_{ij}} using the observation set $\{i^k,A^k\}_{k \in [K(t)]}$, and denote $p_t := \{p_{ij,t}\}_{i \in A_j, j \in [m]} \in \bbR^N$. With today's data collection capabilities, accumulating data in this way and updating the choice frequencies $p_t$ is very realistic. In the limit, they will statistically converge $p_t \to p$ for some $p \in \bbR^N$. Note that if the observations are generated via some `true model' $\lambda^*$, then $p \in X$ (but  our methods will not require this). In the dynamic setting, our goal is still to select $x = \sum_{\sigma \in S_n} \lambda(\sigma) a(\sigma)$ close to $p$, but only with access to the sequence $\{p_t\}_{t \geq 1}$. More precisely, we solve
\vspace{-17.5pt}
\begin{equation}\label{eqn:distance-learn-choice-JEO}
\min_x \left\{ D(x,p):~ x \in X \right\}, \quad \text{given only the sequence } \{p_t\}_{t \geq 1} \text{ s.t. } p_t \to p.
\end{equation}

\subsection{Non-Parametric Choice Model Estimation as a Dynamic Saddle Point Problem}\label{sec:choice-model-saddlept}

In this paper, we consider the following class of general class of distance measures admitting a max-type representation, which casts \eqref{eqn:general-distance-learn-choice} as a min-max saddle point (SP) problem:
\begin{subequations}\label{eqn:SP-choice}
\vspace{-8.5pt}
\begin{align}
D(x,p) &:= \max_{y \in Y} \Psi(x,y;p), \quad \Psi(x,y;p) := \la B(x-p), y \ra - \alpha \omega(y)\label{eqn:distance-max-rep}\\
\SV(p) &:= \min_{x \in X} \max_{y \in Y} \Psi(x,y;p) = \min_{x \in X} D(x,p).\label{eqn:choice-SP}
\end{align}
\end{subequations}
Here, $\alpha \geq 0$, $Y \subseteq \bbR^k$ is closed and convex, $\omega:Y \to \bbR$ is a convex function, and $B$ is a matrix with $N$ columns (the same dimension as $X$) and $k$ rows (the same dimension as $Y$).

\begin{remark}\label{rem:linear-opt-X}
Given any proper convex (in $x$) distance measure $D(x,p)$, through convex conjugacy, a max-type representation \eqref{eqn:distance-max-rep} always exists. We have (for some domain $Y$)
\vspace{-5pt}
\[ D(x,p) = \max_{y \in Y} \left\{ \la y, x \ra - D^*(y, p) \right\}, \quad D^*(y, p) := \max_x \left\{ \la y, x \ra - D(x, p) \right\}. \]
The domain $Y$ is bounded when $\grad_x D(x,p)$ is bounded, and $\min_{x \in X} D(x, p) = \min_{x \in X} \max_{y \in Y} \Psi(x,y; p)$, where $\Psi(x,y; p) := \la y,x \ra - D^*(y, p)$.

Solving $\min_{x \in X} \Psi(x,y;p)$ will be a key component in the algorithms of Sections \ref{sec:JEO} and \ref{sec:JEO-algs}. Our particular form of $\Psi$ in \eqref{eqn:distance-max-rep}, as well as the one that arises from convex conjugacy, indicates that when all else (i.e., $y,p$) is fixed, we are solving a linear optimization problem over $X$. Since $X$ only admits a high-dimensional representation \eqref{eqn:X-choice-domain}, non-linear optimization over $X$ presents difficulties. On the other hand, \emph{linear} optimization over $X$, while non-trivial, is a manageable problem; we elaborate on this subproblem further in Section \ref{sec:comb-subproblem}.

We specifically consider distance measures of the form \eqref{eqn:distance-max-rep} as this form is broad enough to capture a large variety of commonly used $D(x,p)$ %
from practice; see e.g., Table \ref{tab:summary-distances}.\!\!
\epr
\end{remark}

We now discuss desirable properties of $D(x,p)$ and appropriate choices of $Y$, $\omega$ and $B$ in \eqref{eqn:distance-max-rep}. First, note that clearly we have $\Psi(x,y;p)$ is convex in $x$ and concave in $y$ for any choices. However, to ensure that the solution to \eqref{eqn:general-distance-learn-choice} is, in some sense, a `reasonable' solution, we desire the property that, if $p \in X$, then $x = p$ is an optimal solution to \eqref{eqn:general-distance-learn-choice}. An easy  sufficient condition for this is to ensure that $D(x,p) = D(p,x)$, which we state next.
\begin{proposition}\label{prop:symmetric-distance}
	Suppose that the function $h(z) := \max_{y \in Y} \left\{ \la Bz, y \ra - \alpha \omega(y) \right\}$ satisfies $h(z) = h(-z)$ for any vector $z \in \bbR^N$. Then, $D(x,p) = D(p,x) = h(x-p)$, and if in addition $p \in X$, then the optimal solution to \eqref{eqn:general-distance-learn-choice} is $x=p$.
\end{proposition}

On the other hand, symmetry is not a necessary condition to have this reasonable property. We present here one such example which exploits the structure of $X$ and $p$ in the sense that if $x,p \in X$, then $\{x_{ij}\}_{i \in A_j}, \{p_{ij}\}_{i \in A_j}$ live in the simplex $\Delta_{|A_j|}$ for each $j \in [m]$.
\begin{proposition}\label{prop:max-diff-dist}
	Let $B = I$, $Y = \Delta_{|A_1|} \times \ldots \times \Delta_{|A_m|}$, and $\omega(y) = 0$ for all $y \in Y$. Then whenever $p \in X$, $x=p$ solves \eqref{eqn:general-distance-learn-choice}.
\end{proposition}

Even under these assumptions, this form of $D(x,p)$ still captures a large variety of distance measures. We give some examples in Table \ref{tab:summary-distances}. In particular, all norms are covered.\!\!
\begin{table}[t!b!]
	\centering
	\small
	\begin{tabular}{ c | c | c | c }
		$Y$ 		  & $B$ & $\omega(y)$ & $D(x,p)$ \\ \hline
		$\{\|y\|_* \leq 1\} \subset \bbR^N$ & $I_N$ & $0$ & $\|x-p\|$ \\
		$\Delta_{|A_1|} \times \ldots \times \Delta_{|A_m|} \subset \bbR^N$ & $I_N$ & $0$ & $\sum\limits_{j \in [m]} \max\limits_{i \in A_j} (x_{ij} - p_{ij})$ \\
		$\{\|y\|_2 \leq 1\} \subset \bbR^N$ & $I_N$ & $\frac{1}{2} \|y\|^2$ & $\begin{cases}
		\frac{1}{2\alpha}\|x-p\|_2^2, &\text{if } \|x-p\|_2 \leq \alpha\\
		\|x-p\|_2-\frac{\alpha}{2}, &\text{o.w. }\|x-p\|_2 > \alpha
		\end{cases}$ \\
		$\Delta_{2N} \subset \bbR^{2N}$ & $\begin{bmatrix} I_N & -I_N \end{bmatrix}$ & $\sum\limits_{k \in [2N]} y_k \log(y_k)$ & $\alpha \log\left( \sum\limits_{j \in [m]} \sum\limits_{i \in \CA_j} 2\cosh\left( \frac{x_{ij}-p_{ij}}{\alpha} \right) \right)$
	\end{tabular}
	\caption{Examples of different distance measures.}\label{tab:summary-distances}
	\vspace{-5pt}
\end{table}

%% file: DynamicLearning.tex
\section{A Primal-Dual Framework for Dynamic Saddle Point Problems}\label{sec:JEO}

We now build a framework for solving general problems of the form \eqref{eqn:SP-choice}
\vspace{-5pt}
\begin{equation}\label{eqn:general-SP}
\SV(p) := \min_{x \in X} \max_{y \in Y} \Psi(x,y;p),
\end{equation}
where $X,Y$ are convex compact sets, and for every $p$ in the parameter space, $\Psi(x,y;p)$ is convex in $x$ and concave in $y$. When $p$ is fixed, \eqref{eqn:general-SP}
leads to a primal-dual pair of problems:\!
\vspace{-5pt}
\begin{subequations}\label{eq:general-primaldual}
\begin{align}
\CP(p) &:= \min_{x \in X} f(x; p), \quad f(x;p) := \max_{y \in Y} \Psi(x,y; p)\label{eqn:general-primal}\\
\CD(p) &:= \max_{y \in Y} g(y; p), \quad g(y;p) := \min_{x \in X} \Psi(x,y; p)\label{eqn:general-dual}.
\end{align}
\end{subequations}
Note that $f$ and $g$ are convex in $x$ and concave in $y$ respectively, hence \eqref{eqn:general-primal} and \eqref{eqn:general-dual} are convex problems, and by Sion's minimax theorem \citep{Sion1958}, strong duality holds, i.e., $\CP(p) = \CD(p) = \SV(p)$. Note that $f(x;p)$ is exactly $D(x,p)$ from \eqref{eqn:SP-choice}.

Since our main focus is the dynamic data setting where $p$ is unknown but approximated, in \eqref{eqn:general-SP} we make the dependence on the problem data $p$ explicit. More precisely, our goal is to solve \eqref{eqn:general-SP} with access to only a sequence $\{p_t\}_{t \geq 1}$ such that $p^t \to p$, i.e.,
\vspace{-7.5pt}
\begin{equation}\label{eqn:general-SP-JEO}
\SV(p) = \min_{x \in X} \max_{y \in Y} \Psi(x,y;p), \quad \text{given only a sequence } \{p_t\}_{t \geq 1} \text{ s.t. } p_t \to p.
\end{equation}
This is a variant of a joint estimation and optimization (JEO) problem, previously studied in \citet{AhmadiShanbhag2014,HoNguyenKK2016}. However, the previous literature on JEO has studied only the primal problem, rather than the more general SP representation of \eqref{eqn:general-SP}. As we will see, by examining \eqref{eqn:general-SP}, we can devise algorithms with favorable properties for certain settings, specifically addressing the high dimensionality challenge of $X$ in non-parametric choice modeling.

The data sequence $\{p_t\}_{t \geq 1}$ induces a sequence of SP problems $\{\SV(p_t)\}_{t \geq 1}$ with the same structure as \eqref{eqn:general-SP}. In order to develop our framework, one of the main tasks is to relate the sequence back to the original SP problem $\SV(p)$.

To start with, let us consider the traditional (static data) setup when $p$ is known. In this setting, solving \eqref{eqn:general-SP} entails finding a primal-dual pair $(\bar{x},\bar{y}) \in X \times Y$ such that $f(\bar{x}; p) = \CP(p) = \CD(p) = g(\bar{y};p)$, i.e., $\bar{x}, \bar{y}$ are optimal solutions to \eqref{eqn:general-primal},\eqref{eqn:general-dual}, respectively. A key quantity to certify this is the \emph{SP gap}, i.e., the sum of the optimality gaps for \eqref{eqn:general-primal} and \eqref{eqn:general-dual}:
\begin{align}
\epsilonsad^\Psi\left( \bar{x}, \bar{y}; p \right) &:= f(\bar{x}; p) - \CP(p) + \CD(p) - g(\bar{y};p)\notag\\
&= \max_{y \in Y} \Psi(\bar{x},y;p) - \min_{x \in X} \Psi(x,\bar{y};p).\label{eqn:SP-gap}
\end{align}
Note that $\epsilonsad^\Psi\left( \bar{x}, \bar{y}; p \right) \geq 0$ (which is equivalent to weak duality), and $\bar{x},\bar{y}$ are optimal when $\epsilonsad^\Psi\left( \bar{x}, \bar{y}; p \right) = 0$. Moreover, a primal-dual pair $(\bar{x},\bar{y}) \in X \times Y$ satisfying 
$ \epsilonsad^\Psi\left( \bar{x}, \bar{y}; p \right) \leq \epsilon$ 
is referred to be an \emph{$\epsilon$-approximate solution} to the SP problem \eqref{eqn:general-SP}, and then
$\bar{x},\bar{y}$ are \emph{$\epsilon$-optimal solutions} to the primal and dual problems \eqref{eqn:general-primal}, \eqref{eqn:general-dual}, respectively as well.

Our goal, therefore, is to find a primal-dual pair $(\bar{x},\bar{y})$ with certifiably small SP gap. In the JEO setting, however, we must do this without access to the data $p$, but instead use only the sequence $\{p_t\}_{t \geq 1}$. The general framework we employ to do this is as follows:
\begin{enumerate}
\item Fix some time horizon $T \in \bbN$.
\item For each $t \in [T]$:
\begin{itemize}
\item Use $p_t$ to generate a primal-dual pair $(x_t,y_t) \in X \times Y$.
\end{itemize}
\item Select non-negative weights $\{\theta_t\}_{t \in [T]}$ to aggregate the primal-dual sequence $\{x_t,y_t\}_{t \in [T]}$ and obtain a primal-dual pair $(\bar{x}_T^\theta,\bar{y}_T^\theta) \in X \times Y$ as follows:
\vspace{-5pt}
\begin{equation}\label{eqn:primaldual-aggregate}
\Theta_T := \sum_{t \in [T]} \theta_t, \quad \bar{x}_T^\theta := \frac{1}{\Theta_T} \sum_{t \in [T]} \theta_t x_t, \quad \bar{y}_T^\theta := \frac{1}{\Theta_T} \sum_{t \in [T]} \theta_t y_t.
\end{equation}
\item Certify the SP gap $\epsilonsad^\Psi\left( \bar{x}_T^\theta, \bar{y}_T^\theta; p \right)$.
\end{enumerate}

In our framework we deliberately leave steps 2 and 4 vague. To specify these, we present the following fundamental result which relates the data sequence $\{p_t\}_{t \in [T]}$, a generic primal-dual sequence $\{x_t,y_t\}_{t \in [T]}$, and the SP gap $\epsilonsad^\Psi\left( \cdot, \cdot; p \right)$.
\begin{theorem}\label{thm:JEO-SPgap-bound}
For given $T \geq 1$, primal-dual sequence $\{x_t,y_t\}_{t \in [T]}$ and weight sequence $\{\theta_t\}_{t \in [T]}$, let $\bar{x}_T^\theta,\bar{y}_T^\theta$ be defined as in \eqref{eqn:primaldual-aggregate}. Then, given any data sequence $\{p_t\}_{t \in [T]}$, we have
\begin{center}
$ \epsilonsad^\Psi\left( \bar{x}_T^\theta, \bar{y}_T^\theta; p \right) \leq \widehat{\epsilon}\left( \left\{ x_t,y_t; p_t,\theta_t \right\}_{t \in [T]} \right) + \epsilon^\circ\left( \left\{ x_t; p_t,\theta_t \right\}_{t \in [T]} \right) + \epsilon^\bullet\left( \left\{ y_t; p_t,\theta_t \right\}_{t \in [T]} \right),$
\end{center}
\vspace{-5pt}
\begin{subequations}\label{eqn:error-terms}
\begin{align}
\text{where}~~\widehat{\epsilon}\left( \left\{ x_t,y_t; p_t,\theta_t \right\}_{t \in [T]} \right)\!\! &:= \max_{y \in Y} \frac{1}{\Theta_T} \sum_{t \in [T]} \theta_t \Psi(x_t, y; p_t) - \min_{x \in X} \frac{1}{\Theta_T} \sum_{t \in [T]} \theta_t \Psi(x, y_t; p_t)\label{eqn:approx-SP-gap}\\
\epsilon^\circ\left( \left\{ x_t; p_t,\theta_t \right\}_{t \in [T]}; p \right)\!\! &:= \max_{y \in Y} \frac{1}{\Theta_T} \sum_{t \in [T]} \theta_t \left[ \Psi(x_t, y; p) - \Psi(x_t,y; p_t) \right]\label{eqn:primal-error}\\
\epsilon^\bullet\left( \left\{ y_t; p_t,\theta_t \right\}_{t \in [T]}; p \right)\!\! &:= \max_{x \in X} \frac{1}{\Theta_T} \sum_{t \in [T]} \theta_t \left[ \Psi(x, y_t; p_t) - \Psi(x, y_t; p) \right].\label{eqn:dual-error}
\end{align}
\end{subequations}
\end{theorem}
Theorem \ref{thm:JEO-SPgap-bound} states that, given any %
sequence $\{x_t,y_t\}_{t \in [T]}$, we can construct a primal-dual pair via \eqref{eqn:primaldual-aggregate} and certify its SP gap $\epsilonsad^\Psi\left( \bar{x}_T^\theta, \bar{y}_T^\theta; p \right)$ on the true data $p$ by bounding the three upper bound terms $\widehat{\epsilon}, \epsilon^\circ, \epsilon^\bullet$. Thus, Theorem \ref{thm:JEO-SPgap-bound} directly addresses step 4 of the template. We design steps 2 and 3 of our framework to ensure that these upper bound terms are small.

Let us comment on the three upper bound terms. Notice that since we have access to $\{p_t\}_{t \in [T]}$ and we generate the sequences $\{x_t,y_t\}_{t \in [T]}$, we will see that $\{x_t,y_t\}_{t \geq 1}$ will be chosen to make $\widehat{\epsilon}$ small. On the other hand, $p$ is present in both the $\epsilon^\circ, \epsilon^\bullet$ terms, and since $p$ is unknown to us, we do not directly control $\epsilon^\circ, \epsilon^\bullet$. However, these terms become small under a mild regularity assumption on the SP function $\Psi$, and when $p_t \to p$ (which is already assumed in the problem statement \eqref{eqn:general-SP-JEO}). For the rest of this section, we focus on the latter two terms; we will focus on $\widehat{\epsilon}$ in Section \ref{sec:JEO-algs}.

Intuitively, we can think of the latter two terms $\epsilon^\circ, \epsilon^\bullet$ as capturing \emph{the price of estimation}, i.e., they are the errors incurred from using inexact estimates $p_t \neq p$. Indeed, if $p$ were known and we had $p_t = p$ for all $t$, then these terms will disappear. We now give a sufficient uniform Lipschitz continuity condition that ensures $\lim_{T \to \infty} \left\{\epsilon^\circ + \epsilon^\bullet \right\} \leq 0$ whenever $p_t \to p$.
\begin{assumption}\label{ass:JEO-continuity}
Fix a norm $\|\cdot\|$ on the data vectors $p$. There exists a constant $L_\Psi < \infty$ such that for any $(x,y) \in X \times Y$ and data $p,p'$, we have $|\Psi(x,y;p) - \Psi(x,y;p')| \leq L_\Psi \|p - p'\|$.
\end{assumption}
\begin{theorem}\label{thm:JEO-error-convergence}
Suppose that Assumption \ref{ass:JEO-continuity} holds, and that $\{\theta_t\}_{t \geq 1}$ is chosen so that $\theta_t \geq 0$ for all $t \geq 1$, $\Theta_T = \sum_{t \in [T]} \theta_t \to \infty$. Then, whenever $p_t \to p$ for any fixed $p$,
\vspace{-5pt}
\[ \limsup_{T \to \infty} \left[ \epsilon^\circ\left( \left\{ x_t; p_t,\theta_t \right\}_{t \in [T]}; p \right) + \epsilon^\bullet\left( \left\{ y_t; p_t,\theta_t \right\}_{t \in [T]}; p \right) \right] \leq \lim_{T \to \infty} \frac{2L_\Psi}{\Theta_T} \sum_{t \in [T]} \theta_t \|p_t - p\| = 0.\]
\end{theorem}

Theorem \ref{thm:JEO-error-convergence} ensures that the error terms $\epsilon^\circ , \epsilon^\bullet$ will be arbitrarily small if we make $T$ large enough. However, the rate at which these terms converge depends on the rate that $\|p_t - p\| \to 0$. This is a quite natural consequence, since in essence it supports the intuition that the overall performance is limited by the quality of the information $\{p_t\}_{t \geq 1}$ available to us; indeed, it is unreasonable to expect that faster rates are possible without assumptions on the dynamics of the sequence $\{p_t\}_{t \geq 1}$ beyond convergence. Different rates of convergence for $\epsilon^\circ,\epsilon^\bullet$ for common choices of $\theta_t$ (i.e., $\theta_t=1~\forall t$ or $\theta_t=t~\forall t$) are given in Appendix \ref{sec:error-rates}.

In order to get convergence in our framework for the non-parametric choice model estimation problem \eqref{eqn:SP-choice}, Assumption \ref{ass:JEO-continuity} must hold for the $\Psi$ defined in \eqref{eqn:distance-max-rep}. This is indeed the case under a compactness assumption on $Y$, which is satisfied for all examples in Table \ref{tab:summary-distances}.
\begin{proposition}\label{prop:Y-bounded-continuity}
	Let $\Psi(x,y;p) = \la B(x-p), y \ra - \alpha \omega(y)$. Suppose that $Y$ is bounded in some (dual) norm, i.e., $\|y\|_* \leq G_Y < \infty$ for some $G_Y\in\R$. Then Assumption \ref{ass:JEO-continuity} holds for $\Psi$ with $L_\Psi = G_Y \|B\|$, where $\|B\|$ is the operator norm induced by $\|\cdot\|$.
\end{proposition}

%% file: Algorithms.tex
\section{Deriving Algorithms from the Primal-Dual Framework}\label{sec:JEO-algs}

In Section \ref{sec:JEO}, we showed that the SP gap $\epsilonsad^\Psi(\bar{x},\bar{y};p)$ certifies the quality of a primal-dual pair $(\bar{x},\bar{y})$ for the SP problem \eqref{eqn:general-SP}. Then, for the JEO problem \eqref{eqn:general-SP-JEO}, Theorem \ref{thm:JEO-SPgap-bound} provides a bound on the SP gap via three auxiliary terms \eqref{eqn:error-terms}. Theorem \ref{thm:JEO-error-convergence} shows that two of these terms, $\epsilon^\circ$ and $\epsilon^\bullet$, go %
to $0$ as long as a Lipschitz continuity-type assumption on data $p$ (Assumption \ref{ass:JEO-continuity}) is satisfied. Thus, to derive algorithms for \eqref{eqn:general-SP-JEO}, we focus on the first $\widehat{\epsilon}$ term.

Intuitively, the $\widehat{\epsilon}$ term can be interpreted as an approximate SP gap term. Indeed, we have
\vspace{-7.5pt}
\[ \widehat{\epsilon}\left( \left\{ x_t,y_t; p_t,\theta_t \right\}_{t \in [T]} \right) \leq \frac{1}{\Theta_T} \sum_{t \in [T]} \theta_t \left[ \max_{y \in Y} \Psi(x_t, y; p_t) - \min_{x \in X} \Psi(x, y_t; p_t) \right], \]
i.e., $\widehat{\epsilon}$ is upper-bounded by the (weighted) average of the SP gaps defined by $\SV(p_t)$ evaluated at the primal-dual pair $(x_t,y_t)$. This means that, to make $\widehat{\epsilon}$ small, at each time $t \in [T]$, we can choose $(x_t,y_t)$ to simply (approximately) solve a SP problem $\SV(p_t)$ given data $p_t$.

While this is a valid method, even approximately solving $\SV(p_t)$ at each iteration to a given accuracy $\epsilon$ can be expensive, and perhaps intractable. A more careful analysis of $\widehat{\epsilon}$, however, reveals more efficient strategies. Notice that $\widehat{\epsilon}%
$ can be written as
\vspace{-5pt}
\begin{subequations}\label{eqn:SP-gap-two-regret}
	\begin{align}
	\widehat{\epsilon}\left( \left\{ x_t,y_t; p_t,\theta_t \right\}_{t \in [T]} \right) &= \cR_x(\{x_t, y_t;p_t,\theta_t\}_{t\in[T]}) + \cR_y(\{x_t, y_t;p_t,\theta_t\}_{t\in[T]}),\label{eqn:eps-hat-decomp}\\
	\text{where} \quad  
	\cR_x(\{x_t, y_t;p_t,\theta_t\}_{t\in[T]})&:= 
	\frac{1}{\Theta_T} \sum_{t \in [T]} \theta_t \Psi(x_t,y_t; p_t) - \min_{x \in X} \frac{1}{\Theta_T} \sum_{t \in [T]} \theta_t \Psi(x, y_t; p_t) \label{eqn:x-regret}\\
	\cR_y(\{x_t, y_t;p_t,\theta_t\}_{t\in[T]}) &:=
	\max_{y \in Y} \frac{1}{\Theta_T} \sum_{t \in [T]} \theta_t \Psi(x_t, y; p_t) - \frac{1}{\Theta_T} \sum_{t \in [T]} \theta_t \Psi(x_t,y_t; p_t).\label{eqn:y-regret} 
	\end{align}
\end{subequations}

These two terms, $\cR_x$ and $\cR_y$ admit a so-called \emph{regret} interpretation from the \emph{online convex optimization} (OCO) domain; see Section~\ref{sec:regret-algs} and Appendix~\ref{sec:oco-intro} for a brief discussion of OCO. OCO immediately provides us various techniques, in the form of  \emph{regret-minimizing} algorithms, to choose the sequences $\{x_t,y_t\}_{t \geq 1}$ that bound these terms in \eqref{eqn:SP-gap-two-regret}. We now provide a strategy for bounding the two regret terms $\cR_x$ and $\cR_y$ while addressing specific high dimensionality challenge of choice model estimation \eqref{eqn:SP-choice}. As we will see in Section \ref{sec:regret-algs}, typically both regret terms $\cR_x$ and $\cR_y$  can be bounded by $\leq O(1/\sqrt{T})$ after $T$ iterations.

\subsection{Primal Oracle Algorithms for Dynamic Saddle Point Problems}\label{sec:primal-oracle}

To get bound certificates on the SP gap \eqref{eqn:SP-gap}, we must choose a primal-dual sequence $\{x_t,y_t\}_{t \in [T]}$ to minimize the two regret terms \eqref{eqn:x-regret} and \eqref{eqn:y-regret}. Specifically, we choose the primal sequence $\{x_t\}_{t \in [T]}$ to minimize the regret term $\cR_x(\{x_t, y_t;p_t,\theta_t\}_{t\in[T]})$, 
and the dual sequence $\{y_t\}_{t \in [T]}$ to minimize the regret term $\cR_y(\{x_t, y_t;p_t,\theta_t\}_{t\in[T]})$. 
This involves selecting two regret-minimizing algorithms from the OCO literature, and applying them respectively to minimize these two regret terms. While there are many possibilities for doing this, as stated in Remark \ref{rem:linear-opt-X}, for non-parametric choice model estimation we focus on a particular type of combination of regret-minimizing algorithms, what we call \emph{primal oracle algorithms}. This class of algorithms in a sense
bypasses the high-dimensional challenge of the primal domain $X$ (where performing the usual projection or proximal operations onto $X$ are difficult)
and focuses on the dual domain $Y$. Concretely, at each time $t \in [T]$, we will make a lookahead decision $x_t$ \emph{after} choosing $y_t$ by computing
\vspace{-5pt}
\begin{equation}\label{eqn:x-update-lookahead}
x_t \in \argmin_{x \in X} \Psi(x,y_t;p_t).
\end{equation}
This immediately results in the following regret bound.
\begin{lemma}\label{lemma:primal-regret-lookahead}
Let $\{y_t\}_{t \in [T]}$, $\{p_t\}_{t \in [T]}$ and $\{\theta_t\}_{t \in [T]}$ be arbitrary dual, data and weight sequences. Suppose that for each $t \in [T]$, $x_t$ is obtained via \eqref{eqn:x-update-lookahead}. Then, we have
\[ \cR_x(\{x_t, y_t;p_t,\theta_t\}_{t\in[T]}) =  \frac{1}{\Theta_T} \sum_{t \in [T]} \theta_t \Psi(x_t,y_t; p_t) - \min_{x \in X} \frac{1}{\Theta_T} \sum_{t \in [T]} \theta_t \Psi(x, y_t; p_t) \leq 0. \]
\end{lemma}
In fact, with this choice and consequently using Theorem~\ref{thm:JEO-SPgap-bound} and Lemma~\ref{lemma:primal-regret-lookahead}, we deduce that the primal optimality gap is bounded solely by the dual regret term $\cR_y$.
\begin{corollary}\label{cor:primal-JEO-bound}
Let $\{y_t\}_{t \in [T]}$, $\{p_t\}_{t \in [T]}$ and $\{\theta_t\}_{t \in [T]}$ be arbitrary dual, data and weight sequences. Suppose that for each $t \in [T]$, $x_t \in \argmin_{x \in X} \Psi(x,y_t;p_t)$. Denote $\bar{x}_T^\theta := \frac{1}{\Theta_T} \sum_{t \in [T]} \theta_t x_t$. Then, we have
\vspace{-5pt}
\begin{align*}
f(\bar{x}_T^\theta;p) - \min_{x \in X} f(x;p) &\leq \cR_y(\{x_t, y_t;p_t,\theta_t\}_{t\in[T]})
+ \epsilon^\circ\left( \left\{ x_t; p_t,\theta_t \right\}_{t \in [T]}; p \right) + \epsilon^\bullet\left( \left\{ y_t; p_t,\theta_t \right\}_{t \in [T]}; p \right).
\end{align*}
\end{corollary}

\begin{remark}\label{rem:sparsity-guarantee}
When using a primal oracle algorithm \eqref{eqn:x-update-lookahead}, Corollary \ref{cor:primal-JEO-bound} gives a bound on the optimality gap. In the context of choice model estimation, when we use a $D(x,p)$ %
of form \eqref{eqn:distance-max-rep} and solve \eqref{eqn:choice-SP}, solving \eqref{eqn:x-update-lookahead} boils down to a linear program over $X$. This has two important implications. Since $X = \Conv\left( \{a(\sigma)\}_{\sigma \in S_n} \right)$, the optimal solution to \eqref{eqn:x-update-lookahead} is some vertex $a(\sigma)$. Thus, $\bar{x}_T^\theta$ is \emph{readily} expressed as a convex combination of the $\{a(\sigma)\}_{\sigma \in S_n}$. This allows us to easily obtain a distribution $\lambda_T$ over $S_n$ to describe our estimated choice model at iteration $T$.
Furthermore, the overall number of iterations $T$ upper-bounds the sparsity of the learned choice model because each iteration increases the support of $\lambda_T$ by at most $1$ and the initial point can be taken to be any ranking.
As a result, our framework allows us to obtain an explicit dependence between the accuracy parameter $\epsilon$ and the sparsity of the estimated $\epsilon$-optimal model via bounds on the regret term $\cR_y$ and the error terms $\epsilon^\circ, \epsilon^\bullet$.
\epr
\end{remark}

We next discuss some possible regret-minimizing algorithms for $\CR_y$ under two different assumptions on $\Psi$. We also elaborate on the connection of these algorithms to well-known optimization algorithms in the static setup.

\subsection{Regret-Minimizing Algorithms}\label{sec:regret-algs}

In the standard online convex optimization setting, we are given a convex domain $Z$ and a finite time horizon $T$. In each time period $t \in [T]$, the following takes place:
\begin{itemize}
	\item we make a decision $z_t \in Z$ based on \emph{past} information from time steps $1,\ldots,t-1$ only.
	\item Then, a convex loss function $h_t:Z \to \bbR$ is revealed, we suffer loss $h_t(z_t)$ and get some feedback typically in the form of first-order information $\grad h_t(z_t)$.
\end{itemize}
It is usually assumed that the functions $h_t$ are chosen possibly by an all-powerful adversary. As such, our sequence of decisions $\{z_t\}_{t \in [T]}$ is evaluated against the best fixed decision in hindsight, and the (weighted average) difference is defined to be the \emph{weighted regret}:
\vspace{-5pt}
\begin{equation}\label{eqn:weighted-regret}
\frac{1}{\Theta_T} \sum_{t \in [T]} \theta_t h_t(z_t) - \min_{z \in Z} \frac{1}{\Theta_T} \sum_{t \in [T]} \theta_t h_t(z),
\end{equation}
We give a full overview of OCO in Appendix \ref{sec:oco-intro}.

By defining $Z = Y$ and $h_t(\cdot) = -\Psi(x_t,\cdot;p_t)$, the regret term $\cR_y$ in \eqref{eqn:y-regret} becomes \eqref{eqn:weighted-regret}. Thus, we state our algorithms in terms of $h_t$ and $Z$. Throughout this subsection, we fix some norm $\|\cdot\|$ on the Euclidean space in which $Z$ lives, and denote its dual norm by $\|\cdot\|_*$.

\begin{assumption}\label{ass:bounded-gradients}
There exists a $G < \infty$ such that for any $z \in Z$, $\|\grad h(z)\|_* \leq G$.
\end{assumption}

\begin{assumption}\label{ass:prox-setup}
There exists a convex function $\omega:Z \to \bbR$ satisfying: (a) for any $z,z' \in Z$, $\la \grad \omega(z) - \grad \omega(z'), z - z' \ra \geq \|z - z'\|^2$; and (b) there exists $\Omega > 0$ such that for any $z,z' \in Z$, $\omega(z) - \omega(z') \leq \Omega$. %
\end{assumption}

Our first algorithm, \emph{online Mirror Descent (oMD)}, applies to the setting where $h_t$ only satisfies convexity; in particular, it can be non-smooth and non-strongly convex.
\begin{theorem}\label{thm:OMD}
Suppose that $Z$ satisfies Assumption \ref{ass:prox-setup}(a,b), and that for all $t \in [T]$, $h_t$ is convex and that Assumption \ref{ass:bounded-gradients} holds for $h_t$. Given the weight sequence $\{\theta_t\}_{t\in[T]}$, compute
\vspace{-5pt}
\[ z_1 \in \argmin_{z \in Z} \omega(z), \quad z_{t+1} = \argmin_{z \in Z} \left\{ \la \gamma \theta_t \grad h_t(z_t) - \grad \omega(z_t), z \ra + \omega(z) \right\}, \]
where $\gamma = \sqrt{\frac{2 \Omega}{G^2 \sum_{t \in [T]} \theta_t^2}}$. Then
\vspace{-10pt}
\[ \frac{1}{\Theta_T} \sum_{t \in [T]} \theta_t h_t(z_t) - \min_{z \in Z} \frac{1}{\Theta_T} \sum_{t \in [T]} \theta_t h_t(z) \leq \sqrt{2 \Omega G^2 \sum_{t \in [T]} \left(\frac{\theta_t}{\Theta_T}\right)^2}. \]
Furthermore, when $\theta_t = 1$, the upper bound on regret becomes $\sqrt{2 \Omega G^2/T}$.
\end{theorem}

Our next two algorithms work under a relative strong convexity assumption on $h_t$ with respect to $\omega$ that allows us to get improved upper bound on regret.
\begin{assumption}\label{ass:strongly-convex}
Let $\omega:Z \to \bbR$ be a function satisfying Assumption \ref{ass:prox-setup}(a). There exists $\alpha > 0$ such that the function $h - \alpha \omega$ is convex.
\end{assumption}
Under Assumption~\ref{ass:strongly-convex}, oMD algorithm with a different step size and weight sequence results in the following better regret guarantee. 
\begin{theorem}\label{thm:s-OMD}
Suppose that $Z$ satisfies Assumption \ref{ass:prox-setup}(a),
and that for all $t \in [T]$, Assumptions \ref{ass:bounded-gradients} and \ref{ass:strongly-convex} hold for $h_t$. Let $\theta_t = t$ for all $t\in[T]$, and compute
\vspace{-5pt}
\[ z_1 \in \argmin_{z \in Z} \omega(z), \quad z_{t+1} = \argmin_{z \in Z} \left\{ \left\la \frac{2 }{\alpha (\theta_t + 1)} \grad h_t(z_t) - \grad \omega(z_t), z \right\ra + \omega(z) \right\}. \]
Then 
$ \frac{1}{\Theta_T} \sum_{t \in [T]} \theta_t h_t(z_t) - \min_{z \in Z} \frac{1}{\Theta_T} \sum_{t \in [T]} \theta_t h_t(z) \leq \frac{2 G^2}{\alpha(T+1)}. $
\end{theorem}

Let us also examine the \emph{follow-the-leader (FTL)} algorithm from OCO literature.
\begin{theorem}\label{thm:s-FTL}
Suppose that $Z$ satisfies Assumption \ref{ass:prox-setup}(a),
and that for all $t \in [T]$, Assumption \ref{ass:strongly-convex} holds for $h_t$, and Assumption \ref{ass:bounded-gradients} holds for $h_t - \alpha \omega$. Let $\theta_t = t$ for all $t\in[T]$, and compute
\vspace{-13pt}
\[ z_1 \in \argmin_{z \in Z}\, \omega(z), \quad z_{t+1} = \argmin_{z \in Z} \left\{ {1\over \Theta_t} \sum_{s \in [t]} \theta_s h_s(z) \right\}. \]
Then 
$ \frac{1}{\Theta_T} \sum_{t \in [T]} \theta_t h_t(z_t) - \min_{z \in Z} \frac{1}{\Theta_T} \sum_{t \in [T]} \theta_t h_t(z) \leq \frac{2 G^2}{\alpha(T+1)}. $
\end{theorem}

For common domains $Z$ such as simplex, Euclidean ball, and spectahedron, standard selections of norm $\|\cdot\|$, distance generating function $\omega(\cdot)$, and set width $\Omega$ satisfying Assumption~\ref{ass:prox-setup}(a,b) and the resulting update rule for $z_{t+1}$ computation used in  %
Theorems~\ref{thm:OMD} and \ref{thm:s-OMD} are discussed by \citet[Section 5.7]{JuditNem2012Pt1}.

Finally, we give conditions on \eqref{eqn:distance-max-rep}, \eqref{eqn:choice-SP} to ensure that Theorems \ref{thm:OMD}, \ref{thm:s-OMD} or \ref{thm:s-FTL} can be applied to minimize the dual regret \eqref{eqn:y-regret}. Note that for our class of distance measures \eqref{eqn:distance-max-rep}, because by definition $X$ and $p$ are bounded (see \eqref{eqn:pq_{ij}} and \eqref{eqn:X-choice-domain}), Assumptions \ref{ass:bounded-gradients} and \ref{ass:strongly-convex} are satisfied depending solely on $\omega$. Thus, a judicious definition of $\omega$ will allow us to utilize the regret minimizing algorithms of Theorems \ref{thm:OMD}, \ref{thm:s-OMD} or \ref{thm:s-FTL}.
\begin{proposition}\label{prop:check-assumptions}
	For any $x \in X$, let $h_x(z) = -\Psi(x,z;p) := \la B(p-x), z \ra + \alpha \omega(z)$. If there exists a norm $\|\cdot\|$ and some $G' > 0$ such that $\|\grad \omega(z)\|_* \leq G'$ for all $z \in Y$, then Assumption \ref{ass:bounded-gradients} holds for $h_x$. If $\omega$ satisfies Assumption \ref{ass:prox-setup}(a) and $\alpha > 0$, then Assumption \ref{ass:strongly-convex} holds for $h_x$.
\end{proposition}

\begin{remark}\label{rem:rates}
Proposition \ref{prop:check-assumptions} allows us to use any of the regret-minimizing algorithms from Theorems \ref{thm:OMD}, \ref{thm:s-OMD} or \ref{thm:s-FTL} to bound the dual regret $\CR_y$, from which primal optimality gap bounds can be inferred through \eqref{eqn:SP-gap} (since the dual optimality gap is always non-negative). Ignoring for the moment the error terms $\epsilon^\circ, \epsilon^\bullet$, using Theorem \ref{thm:OMD} results in a point $\bar{x}_T^\theta$ with $\epsilon$-gap after $O(1/\epsilon^2)$ iterations. By Remark \ref{rem:sparsity-guarantee}, the support of the learned choice model is also at most $O(1/\epsilon^2)$, and thus our framework simultaneously results in both errors bounds and \emph{explicit} sparsity guarantees. See Appendix \ref{sec:FWcomparison-rates-MD} and \ref{sec:FWcomparison-rates-FW} for more details on the rates attainable using our framework, and Appendix~\ref{sec:error-rates} for rates for the error terms $\epsilon^\circ, \epsilon^\bullet$.
\epr
\end{remark}

\begin{remark}\label{rem:comparison-FW}
Appendix \ref{sec:FW-nonsmooth} shows that the distances measures $D$ used in the previous literature are indeed non-smooth, which is also the case when $\alpha=0$ in \eqref{eqn:SP-choice}. \citet[Example 1]{Nesterov2018} shows that applying the Frank-Wolfe (F-W) algorithm as suggested by \citet{JagabathulaRusmevichientong2016} to a non-smooth objective does not converge in general. When we have $\alpha > 0$ in \eqref{eqn:distance-max-rep} and $\omega(y)$ is strongly convex (Assumption \ref{ass:prox-setup}(a)), then $D(x,p)$ is smooth in $x$. Also, when $D(x,p) = \|x-p\|_q$ for $q \in [2,\infty)$, then simply squaring the norm makes $D(x,p)$ smooth. This means that the F-W algorithm can be applied for both of these smooth cases. \citet{FreundGrigas2016} provide results for F-W with approximate gradients, and in Appendix \ref{sec:FWcomparison-rates-FW}, we show that in the dynamic setting, $\grad_x D(x,p_t)$ can be interpreted as an approximate gradient for $\grad_x D(x,p)$, thus the results of \citet{FreundGrigas2016} can be applied here. However, we also show that ensuring the convergence of F-W with approximate gradients using their results requires an assumption that $\|p - p_t\| \to 0$ \emph{sufficiently fast}. In contrast, our primal-dual framework bypasses such assumptions on smoothness of $D(x,p)$ and data convergence rate, and works for general $D(x,p)$, e.g., even non-smooth cases of $D(x,p) = \|x-p\|_q$ for $q \in \{1,\infty\}$.
\epr
\end{remark}

\begin{remark}\label{rem:interpretations}
In Appendix \ref{sec:interpretations}, in the static setting ($p_t=p$), we show that our primal oracle algorithms, i.e., ones which use \eqref{eqn:x-update-lookahead} for $x$ in combination with Theorems \ref{thm:OMD}, \ref{thm:s-OMD} or \ref{thm:s-FTL} for $y$, can be interpreted as running traditional MD on a particular dual problem, or F-W on $D(x,p)$. By placing this into our primal-dual framework of Section \ref{sec:JEO}, we are able to obtain guarantees for these algorithms in the dynamic setting as well. Furthermore, in the dynamic setting, the F-W algorithm derived from our framework is slightly different to the traditional one from \citet{FreundGrigas2016} discussed in Remark \ref{rem:comparison-FW} and Section \ref{sec:FWcomparison-rates-FW}. In particular, the analysis of our F-W variant through the primal-dual framework does not require that $\|p - p_t\| \to 0$ at any particular rate.
In Section \ref{sec:computationSummary}, we also numerically compare several algorithms derived from our framework (using Theorems \ref{thm:OMD}, \ref{thm:s-OMD} and \ref{thm:s-FTL}) with ones using the traditional F-W algorithm of \citet{FreundGrigas2016} %
on 
the non-smooth distance function $D(x,p)=\|x-p\|_2$ and an appropriate smoothed version for methods that only work on smooth objectives.
We find that algorithms that work directly with the non-smooth $D(x,p)$ outperform their counterparts requiring smoothness.
\epr
\end{remark}

\subsection{Combinatorial Subproblem}\label{sec:comb-subproblem}

One of the key steps in our framework is solving the LO problem \eqref{eqn:x-update-lookahead} over $X$ (see also the definition of $\Psi$ in \eqref{eqn:choice-SP}). By \eqref{eqn:X-choice-domain}, $X$ is a polytope with vertices $a(\sigma)$, we have, for a cost vector $c$,
\vspace{-8pt}
\[ \argmin_{x \in X} \la c,x \ra = \Conv\left(\left\{ a(\sigma^*) : \sigma^* \in \argmin_{\sigma \in S_n} \la a(\sigma), c \ra \right\}\right). \]
Thus, in each iteration, we must solve the following combinatorial optimization problem over rankings:
\vspace{-8pt}
\begin{equation}\label{eqn:subgradient-subproblem}
\min_{\sigma} \left\{ \sum_{j \in [m]} \sum_{i \in A_j} y_{ij}^t a_{ij}(\sigma) :~ \sigma \in S_n \right\}.
\end{equation}
Problem~\eqref{eqn:subgradient-subproblem} is NP-hard, since it is a generalization of the linear ordering problem and the maximum weighted independent set problem, see e.g., \cite[Proposition 3]{vanRyzinVulcano2015}. However, we note that the exact same combinatorial problem must be solved in all other approaches of learning a non-parametric choice model (see Appendix \ref{sec:learning-choice-model-existing} and in particular the equations \eqref{eqn:revenue-prediction-robust-constraint}, \eqref{eqn:MLE-column-generating}, \eqref{eqn:norm-l1-min-column-generating}). While we cannot avoid the NP-hardness in learning a non-parametric choice model from data, we note that \eqref{eqn:subgradient-subproblem} can be formulated as a (relatively) compact integer program with $O(n^2)$ variables and $O(n^3)$ constraints, and also it can be handled efficiently by off-the-shelf integer programming solvers; see Figure~\ref{fig:metrics-static-fixedthread} in Appendix~\ref{sec:computationApp}.
Furthermore, \citet{JagabathulaRusmevichientong2016} prove polynomial-time solvability of \eqref{eqn:subgradient-subproblem} under a number of assumptions on the subsets $A_1,\ldots,A_m$.

%% file: ExperimentsMain.tex
\section{Computational Study}\label{sec:computationSummary}
We carried out a computational study to examine the performance of various methods for estimating a non-parametric choice model from dynamic data. Full details of our computational setup and empirical observations are given in Appendix \ref{sec:experiments}. For our results in this section, we considered the $\ell_2$-norm based distance measure, i.e., $D(x,p) = \|x-p\|_2$.

We tested four algorithms derived from our framework: a F-W method using the smoothed norm ($\text{FW}_\text{dyn}$), a dual MD method on the natural (non-smooth) norm ($\text{MD}_\text{ns}$), the same MD method on the smoothed norm ($\text{MD}_\text{smth}$) as well as the squared $\ell_2$-norm ($\text{MD}_\text{sq}$). We also tested two other potential algorithms: a na\"{i}ve F-W method on the smoothed norm ($\text{FW}_\text{na\"{i}ve}$), and a (na\"{i}ve) F-W method on the squared $\ell_2$-norm ($\text{FW}_\text{sq}$). See Appendix \ref{sec:experiments} for details on each method.

Observations are generated from a ground truth mixed multinomial logit model. Observations are constrained to belong to some fixed `training' set of item-subset pairs. For each iteration in each of the methods, we generated $\kappa \in \bbN$ new observations and updated $p_t$ (using both old and new observations) according to \eqref{eqn:pq_{ij}}.

Our stopping criterion for all algorithms is the same: the mean absolute error (MAE) of the choice probabilities versus the empirical frequencies on the training item-subset pairs should be below a set threshold (0.001). We evaluate the performance of each method using two metrics: the MAE on a `test' set of unseen item-subset pairs and the number of iterations until termination. We found that the number of different observed rankings at termination (i.e., the model sparsity) is strongly correlated with the number of iterations. 

\begin{figure}[t!bh]
\centering
\includegraphics[page=1,scale=1]{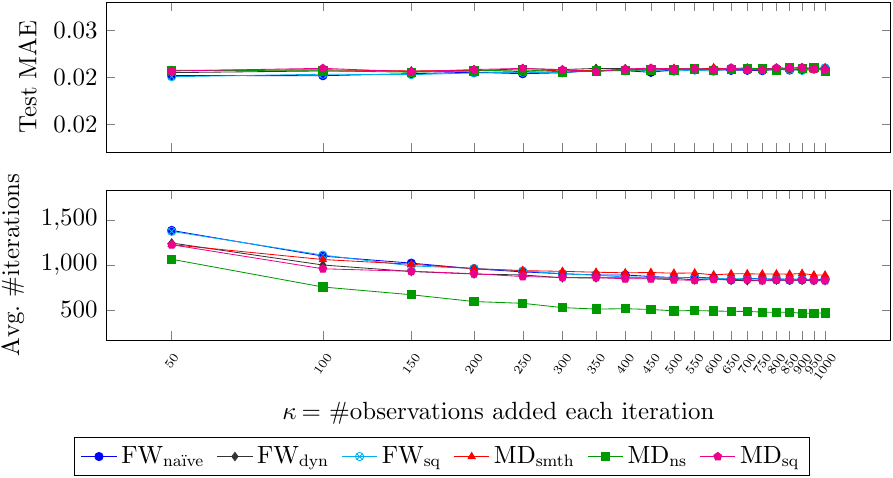}
\caption{Performance metrics using $\ell_2$-norm for different methods  for the ground truth choice model with $n=10$, $m=20$, $K=L=5$.}
\label{fig:metrics-dynamic}
\vspace{-10pt}
\end{figure}

Figure \ref{fig:metrics-dynamic} compares each method tested as well as the effect of using different $\kappa$ on our performance metrics, namely the test MAE and the number of iterations (we note that the solution times are highly correlated with the number of iterations; see Figure \ref{fig:metrics-static-fixedthread} and Figure \ref{fig:compare_m}(b)). In Figure \ref{fig:metrics-dynamic}, different methods and $\kappa$ have very little effect on the final test set MAE of the estimated model. However, the number of iterations until termination decreases as $\kappa$ increases. Furthermore, we see that the $\text{MD}_\text{ns}$ method performs the best in terms of the number of iterations across all levels of $\kappa$. Notice also that, while the difference is minor, the $\text{FW}_\text{dyn}$ method derived from our framework outperforms the $\text{FW}_\text{na\"{i}ve}$ method. Both of these highlight the benefit of using our primal-dual framework to derive optimization algorithms for estimating non-parametric choice models. Theoretically, even though all methods have the same asymptotic $O(1/\sqrt{T})$ suboptimality gap bound, the $\text{MD}_\text{ns}$ method treats the original non-smooth problem directly, and enjoys superior constant factors in the convergence rate. We believe this may be the main factor contributing to its superior numerical performance.

\begin{figure}[t!bh]
\centering
\includegraphics[page=2,scale=1]{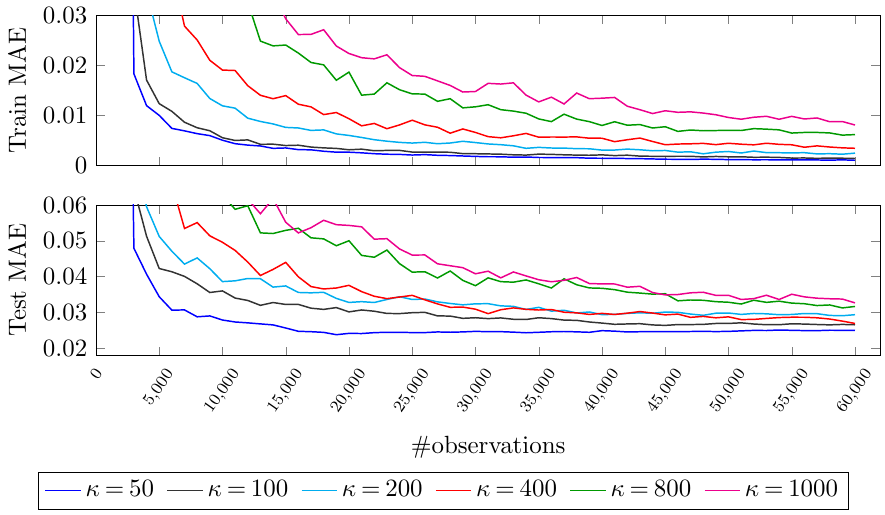}
\caption{Train and test set MAE vs \#observations seen, using $\ell_2$-norm minimization with (non-smooth) dual MD, i.e., $\text{MD}_{\text{ns}}$ algorithm, for the ground truth choice model with $n=10$, $m=20$, $K=L=5$.}
\label{fig:metrics-dynamic-obs}
\vspace{-10pt}
\end{figure}

Figure \ref{fig:metrics-dynamic-obs} examines the training and test set MAE after seeing a certain number of observations, for one particular instance, specifically for the $\text{MD}_\text{ns}$ method. We compare the profiles for different $\kappa$ values to see the impact on MAE of the rate at which new observations are utilized. Note that, at a given number of observations $K\in\{5000,10000,\ldots,60000\}$, if $\kappa$ is higher, then necessarily the method will have performed less iterations. Thus, somewhat unsurprisingly, we observe that both training and test set MAE are lower for lower $\kappa$ for any given number of observations $K$. Figures \ref{fig:metrics-dynamic} and \ref{fig:metrics-dynamic-obs} together suggests the following. Suppose we determine that we perform iterates at consistent time intervals (e.g., two iterations per second). In the regime where observations arrive slowly (e.g., $\kappa = 50$), Figure \ref{fig:metrics-dynamic-obs} suggests that we need not wait longer for more observations to achieve low MAE. In the regime where observations arrive quickly (e.g., $\kappa = 1000$), Figure \ref{fig:metrics-dynamic} suggests that we need not speed up our iteration updates in order to converge quicker.

%% file: Conclusion.tex
\section{Conclusions}\label{sec:conclusion}

In this paper, by studying a general convex-concave SP JEO problem and utilizing OCO, we present an efficient and unified framework to solve non-parametric choice estimation problem with dynamic data. In both the static and dynamic settings, our framework immediately provides error bounds and convergence guarantees (on the number of iterations needed to achieve a certain estimation accuracy), which in turn provide guarantees on the sparsity of our estimated choice model.
Our computational experiments %
demonstrate improved performance when algorithms derived from our framework are used.

In terms of the convex optimization literature, our developments on solving non-smooth optimization problems facing a high-dimensionality challenge and their dynamic variants through JEO for the SP problems are new as well. For static \emph{smooth} problems with a high dimensionality challenge, usually linear optimization as opposed to projection is more tractable, and in such cases the F-W algorithm is often used. In both static and dynamic settings, our framework provides several alternatives to this, all utilizing a linear optimization oracle, but also enjoying convergence in the more general non-smooth setup, which the F-W algorithm does not.

%% file: Proofs.tex
\section{Proofs}\label{sec:proofs}

\proof{Proof of Proposition \ref{prop:symmetric-distance}.}
Note that $D(x,p) = h(x-p) = D(p,x)$. Furthermore, $h(z)$ is clearly convex in $z$, and hence $h(0) \leq h(z)/2 + h(-z)/2 = h(z)$ for any $z$, so $x=p$ is an optimal solution.
\Halmos
\endproof

\proof{Proof of Proposition \ref{prop:max-diff-dist}.}
It is easy to see that
\[ D(x,p) = \sum_{j \in [m]} \max_{i \in A_j} (x_{ij} - p_{ij}). \]
Now, when $x = p$, $D(x,p) = 0$. However, when $x \neq p$, there exists at least one $j \in [m]$ and $i \in A_j$ such that $x_{ij} \neq p_{ij}$. If $x_{ij} > p_{ij}$, then we know $D(x,p) > 0$. If $x_{ij} < p_{ij}$ then since both $\{x_{ij}\}_{i \in A_j}, \{p_{ij}\}_{i \in A_j} \in \Delta_{|A_j|}$, we have $\sum_{i \in A_j} (x_{ij} - p_{ij}) = 0$, so there must exist some other $i' \in A_j$ such that $x_{i'j} > p_{i'j}$, hence $D(x,p) > 0$ also.
\Halmos
\endproof

\proof{Proof of Theorem \ref{thm:JEO-SPgap-bound}.}
First, notice that $\{\theta_t/\Theta_T\}_{t \in [T]}$ form a set of convex combination weights. Thus, using the convex-concave structure of $\Psi(\cdot,\cdot; p)$ and the definition of $\epsilonsad^\Psi\left( \bar{x}_T^\theta, \bar{y}_T^\theta; p \right)$ in 
\eqref{eqn:SP-gap}, %
we get the standard bound
\[\epsilonsad^\Psi\left( \bar{x}_T^\theta, \bar{y}_T^\theta; p \right) \leq \max_{y \in Y} \frac{1}{\Theta_T} \sum_{t \in [T]} \theta_t \Psi(x_t, y; p) - \min_{x \in X} \frac{1}{\Theta_T} \sum_{t \in [T]} \theta_t \Psi(x, y_t; p).\]
Let us examine the first term on the right hand side. Adding and subtracting $\Psi(x_t,y; p_t)$ for each term in the sum, we can bound this term by
\begin{align*}
\max_{y \in Y} \frac{1}{\Theta_T} \sum_{t \in [T]} \theta_t \Psi(x_t, y; p) &= \max_{y \in Y} \frac{1}{\Theta_T} \sum_{t \in [T]} \theta_t \left[ \Psi(x_t, y; p_t) + \Psi(x_t, y; p) - \Psi(x_t, y; p_t) \right]\\
&= \max_{y \in Y} \left\{ \frac{1}{\Theta_T} \sum_{t \in [T]} \theta_t \Psi(x_t, y; p_t) + \frac{1}{\Theta_T} \sum_{t \in [T]} \theta_t \left[ \Psi(x_t, y; p) - \Psi(x_t, y; p_t) \right] \right\}\\
&\leq \max_{y \in Y} \frac{1}{\Theta_T} \sum_{t \in [T]} \theta_t \Psi(x_t, y; p_t) + \max_{y \in Y} \frac{1}{\Theta_T} \sum_{t \in [T]} \theta_t \left[ \Psi(x_t, y; p) - \Psi(x_t, y; p_t) \right]\\
&= \max_{y \in Y} \frac{1}{\Theta_T} \sum_{t \in [T]} \theta_t \Psi(x_t, y; p_t) + \epsilon^\circ\left( \left\{ x_t; p_t,\theta_t \right\}_{t \in [T]}; p \right).
\end{align*}
Using a similar strategy for the second term $\min_{x \in X} \frac{1}{\Theta_T} \sum_{t \in [T]} \theta_t \Psi(x, y_t; p)$, we can get
\[\min_{x \in X} \frac{1}{\Theta_T} \sum_{t \in [T]} \theta_t \Psi(x, y_t; p) \geq \min_{x \in X} \frac{1}{\Theta_T} \sum_{t \in [T]} \theta_t \Psi(x, y_t; p_t) - \epsilon^\bullet\left( \left\{ y_t; p_t,\theta_t \right\}_{t \in [T]}; p \right).\]
Subtracting this lower bound from the upper bound on the first term then gives us the result.
\Halmos
\endproof

\proof{Proof of Theorem \ref{thm:JEO-error-convergence}.}
First observe that for any $t \geq 1$, $(x_t,y), (x,y_t) \in X \times Y$ and $p$, by Assumption \ref{ass:JEO-continuity},
\[ |\Psi(x_t,y; p) - \Psi(x_t,y; p_t)| \leq L_\Psi \|p - p_t\|, \quad |\Psi(x,y_t; p_t) - \Psi(x,y_t; p)| \leq L_\Psi \|p_t - p\|. \]
This implies that
\[ \epsilon^\circ\left( \left\{ x_t; p_t,\theta_t \right\}_{t \in [T]}; p \right) + \epsilon^\bullet\left( \left\{ y_t; p_t,\theta_t \right\}_{t \in [T]}; p \right) \leq \frac{2L_\Psi}{\Theta_T} \sum_{t \in [T]} \theta_t \|p_t - p\|. \]

We now show the following:
\[ a_t \to 0 \implies \frac{1}{\Theta_T} \sum_{t \in [T]} \theta_t a_t \to 0. \]
To get our result, we apply this to the sequence $a_t = 2 L_\Psi \|p_t - p\|$, which converges to $0$ since $p_t \to p$. Fix some $\epsilon > 0$, and choose $S(\epsilon) \in \bbN$ sufficiently large such that for $t \geq S(\epsilon)$, $|a_t| \leq \epsilon/3$. Furthermore, choose $T$ sufficiently large such that $\left| \frac{1}{\Theta_T} \sum_{t \in [S(\epsilon)]} \theta_t a_t \right| \leq \epsilon/2$ and $\left|\frac{1}{\Theta_T} \sum_{t \in [S(\epsilon)]} \theta_t \right| \leq 1/2$. We have
\[\left| \frac{1}{\Theta_T} \sum_{t \in [T]} \theta_t a_t \right| \leq \left| \frac{1}{\Theta_T} \sum_{t \in [S(\epsilon)]} \theta_t a_t \right| + \left| \frac{1}{\Theta_T} \sum_{t=S(\epsilon)+1}^T \theta_t a_t \right|.\]
The first term is $\leq \epsilon/2$ by our choice of $T$, and also the second term satisfies 
\[\left|\frac{1}{\Theta_T} \sum_{t=S(\epsilon)+1}^T \theta_t a_t \right| \leq \frac{1}{\Theta_T} \sum_{t=S(\epsilon)+1}^T \theta_t |a_t| \leq \frac{\epsilon}{3\Theta_T} \sum_{t=S(\epsilon)+1}^T \theta_t = \frac{\epsilon}{3} \left(1 - \frac{1}{\Theta_T} \sum_{t=1}^{S(\epsilon)} \theta_t \right) \leq \frac{\epsilon}{2}.\]
\Halmos
\endproof

\proof{Proof of Proposition \ref{prop:Y-bounded-continuity}.}
We have $|\Psi(x,y;p) - \Psi(x,y;p')| = |\la B(p-p'), y \ra| \leq \|y\|_* \|B (p-p')\| \leq \|y\|_* \|B\| \|p - p'\| \leq G_Y \|B\| \|p - p'\|.$
\Halmos
\endproof

\proof{Proof of Lemma \ref{lemma:primal-regret-lookahead}.}
This is immediate since for any $x \in X$, $\Psi(x_t,y_t; p_t) = \min_{x' \in X} \Psi(x',y_t; p_t) \leq \Psi(x,y_t;p_t)$.
\Halmos
\endproof

\proof{Proof of Corollary \ref{cor:primal-JEO-bound}.}
Theorem \ref{thm:JEO-SPgap-bound} and Corollary \ref{cor:primal-JEO-bound} immediately imply that $\epsilonsad^\Psi(\bar{x}_T^\theta, \bar{y}_T^\theta; p)$ is bounded by the right hand side in the result. To get the left hand side, notice that
\[ \epsilonsad^\Psi(\bar{x}_T^\theta, \bar{y}_T^\theta; p) = f(\bar{x}_T^\theta;p) - \min_{x \in X} f(x;p) + \max_{y \in Y} g(y,p) - g(\bar{y}_T^\theta;p) \geq f(\bar{x}_T^\theta;p) - \min_{x \in X} f(x;p) \]
since the optimality gap for the dual problem $\CD(p)$ in \eqref{eqn:general-dual} is always non-negative.
\Halmos
\endproof

\proof{Proof of Theorems \ref{thm:OMD}, \ref{thm:s-OMD} and \ref{thm:s-FTL}.}
Theorem \ref{thm:OMD} follows almost directly from \citetEC[Theorem 1]{HoNguyenKK2016}, with minor modifications. Theorem \ref{thm:s-OMD} follows from \citetEC[Theorem 2]{HoNguyenKK2016}. Theorem \ref{thm:s-FTL} follows from \citetEC[Theorem 2]{ShalevShwartzKakade2008} which, after taking $\ell_t = \theta_t \alpha \omega + \theta_t (h_t - \alpha \omega)$ and $\theta_t = t$, gives the regret bound
\[ \sum_{t \in [T]} \theta_t h_t(z_t) - \min_{z \in Z} \sum_{t \in [T]} \theta_t h_t(z) \leq \frac{1}{2} \sum_{t \in [T]} \frac{\theta_t^2 G^2}{\alpha \sum_{s \in [t]} \theta_s} = \frac{1}{\alpha} \sum_{t \in [T]} \frac{t G^2}{t+1} \leq \frac{TG^2}{\alpha}. \]
Dividing by $\Theta_T = T(T+1)/2$ gives the result.
\Halmos
\endproof

\proof{Proof of Proposition \ref{prop:check-assumptions}}
Note that $\grad_z h_x(z) = B(p-x) + \alpha \omega(z)$, hence $\|\grad_z h_x(z)\|_* \leq \|B(p-x)\|_* + G'$ is uniformly bounded over $p,x$ since these come from a bounded set, so Assumption \ref{ass:bounded-gradients} is satisfied. Assumption \ref{ass:strongly-convex} holds trivially.
\Halmos
\endproof

%% file: ExistingApproaches.tex
\section{Existing Approaches to Non-Parametric Choice Estimation}\label{sec:learning-choice-model-existing}

In this appendix, we examine the existing approaches to learn the non-parametric choice model, i.e., infer an appropriate probability vector $\lambda$ using the data collected via the process outlined in Section~\ref{sec:model-data}, and demonstrate how they are particular instantiations of our general model. For a fixed subset $\CA_j$, $j \in [m]$, we denote the collection of associated choice probabilities as $A_j \lambda = \{\bbP_{\lambda}[i \mid \CA_j]\}_{i \in \CA_j} \in \Delta_{|\CA_j|}$.

\subsection{Revenue Prediction Approach}\label{sec:rev-predict}

Let $r_i$ be the revenue of item $i\in[n]$. Then the expected revenue of an assortment $\cA \subset [n]$ under distribution $\lambda$ is $\sum_{i \in \cA} r_i \bbP_{\lambda}[i \mid \cA]$.
\citetEC{FariasJS2013} seek to find the worst-case expected revenue from a distribution $\lambda$ consistent with the given data in the sense that the theoretical probabilities $\bbP_\lambda[i \mid \cA_j] = \la a_{ij}, \lambda \ra$ are precisely consistent with their empirical estimates $p_{ij}$. Since the probabilities $\bbP_{\lambda}[i \mid \cA]$ are linear in $\lambda$, this can be formulated as a linear program (LP)
\begin{align*}
\min_{\lambda} \left\{ \sum_{i \in \cA} r_i \bbP_{\lambda}[i \mid \cA] :~ A\lambda = p,~ \lambda \in \Delta_{n!}\right\}.
\end{align*}
We first make a few observations related to this model of \citetEC{FariasJS2013}. 
In fact, when $\cA = \cA_j$ for some $j \in [m]$, we have $\bbP_{\lambda}[i \mid \cA] = \la a_{ij}, \lambda \ra = p_{ij}$ due to the constraints $A\lambda = p$, hence the objective is constant. Thus, the LP becomes a feasibility problem
\begin{align}\label{eqn:revenue-prediction-model}
\text{find} \quad \lambda \in \Delta_{n!} \quad \text{s.t.} \quad  A\lambda = p.
\end{align}
That said, \eqref{eqn:revenue-prediction-model} is still computationally intractable even for moderate values of $n$ because it involves $n!$ variables. Nonetheless, the dual of \eqref{eqn:revenue-prediction-model} admits the following robust LP interpretation:
\begin{align}\label{eqn:revenue-prediction-dual}
\max_{\beta,\nu} \left\{ \la \beta, p\ra - \nu :~ \max_{\sigma \in S_n} \la \beta, a(\sigma)\ra \leq \nu \right\}.
\end{align}
Note that verifying the feasibility of a solution with respect to the robust constraint in \eqref{eqn:revenue-prediction-dual}, i.e., 
\begin{equation}\label{eqn:revenue-prediction-robust-constraint}
\max_{\sigma \in S_n} \la \beta, a(\sigma) \ra= \max_{\sigma} \left\{ \sum_{j \in [m]} \sum_{i \in \cA_j} \beta_{ij} a_{ij}(\sigma) :~ \sigma \in S_n \right\} \leq \nu
\end{equation}
is a combinatorial problem of the exact same form as \eqref{eqn:subgradient-subproblem}.
\citetEC{FariasJS2013} suggests solving \eqref{eqn:revenue-prediction-dual} either using the constraint sampling technique \citepEC{CalafioreCampi2005} or by building an approximation to its robust counterpart obtained from approximating the uncertainty sets with an efficiently representable polyhedron.

In fact, \eqref{eqn:revenue-prediction-model} can be seen as choosing $\lambda \in \Delta_{n!}$ to minimize a (very harsh) distance measure:
\begin{equation}\label{eqn:revenue-prediction-distance-min}
\min_{\lambda \in \Delta_{n!}} D(A\lambda,p), \quad D(A\lambda,p) = \begin{cases} 0, & A\lambda = p\\ \infty, & \text{otherwise}. \end{cases}
\end{equation}
In general, and specifically when the observations are noisy, there is no guarantee that there exists $\lambda \in \Delta_{n!}$ to fit the data $p$ exactly, i.e., $A\lambda = p$. To remedy this, \citetEC{vanRyzinVulcano2015} and \citetEC{BertsimasMisic2015} examine approaches that use less harsh distance measures $D(\cdot,\cdot)$.

\subsection{Maximum Likelihood Estimation Approach}\label{sec:MLE-approach}

\citetEC{vanRyzinVulcano2015} propose the following method to learn $\lambda$ via maximum likelihood estimation (MLE). We next describe their method and provide an alternative interpretation of their approach as the minimization of a particular distance measure, namely  Kullback-Leibler (KL) divergence, between the true distributions $A_j \lambda$ and their empirical estimates $p_j$. 
Note that given two positive vectors $p,x\in\R^n$, their KL divergence is $\KL(p,x):=\sum_{i\in[n]}p_i\log(x_i/p_i)$.

By \eqref{eqn:pq_{ij}}, each item-assortment pair $i \in \CA_j$ is seen $Kq_{ij}$ times amongst the observations $\left\{i^k,\cA^k\right\}_{k=1}^K$. Based on this, the log-likelihood of the observation set $\left\{i^k,\cA^k\right\}_{k=1}^K$ is
$\sum_{j \in [m]} \sum_{i \in \cA_j} Kq_{ij} \log\left( \la a_{ij}, \lambda \ra\right)$. 
Thus, ignoring the constant $K$ factor, the MLE problem is
\begin{align}\label{eqn:MLE-frequency}
\max_\lambda \left\{ \sum_{j \in [m]} \sum_{i \in \cA_j} q_{ij} \log\left( \la a_{ij}, \lambda \ra\right) :~  \lambda \in \Delta_{n!} \right\}.
\end{align}
Throughout, we use the convention that when $q_{ij} = \la a_{ij}, \lambda \ra= 0$, we set $q_{ij} \log(\la a_{ij}, \lambda\ra) = 0$. This implies that if the optimal solution $\lambda$ to \eqref{eqn:MLE-frequency} has $\bbP_\lambda[i \mid \cA_j] = \la a_{ij}, \lambda \ra= 0$, then we must have $q_{ij} = 0$ also, i.e., we did not observe any choices of $i$ from $\cA_j$ in our data either.

Like \eqref{eqn:revenue-prediction-model}, the problem \eqref{eqn:MLE-frequency} is very large, with $n!$ variables. A column generation technique is suggested
in \citetEC{vanRyzinVulcano2015} 
to get around this, i.e., solve \eqref{eqn:MLE-frequency} on a subset of the variables, and use the optimality conditions to add variables as needed.
The MLE column generating subproblem is constructed as
\begin{equation}\label{eqn:MLE-column-generating}
\max_{\sigma} \left\{ \sum_{j \in [m]} \sum_{i \in \cA_j} \frac{q_{ij} a_{ij}(\sigma)}{\la a_{ij}, \lambda(S)\ra} :~ \sigma \in S_n \right\}.
\end{equation}
The solution $\lambda(S)$ is optimal if $\eqref{eqn:MLE-column-generating} \leq K$, otherwise the column $\sigma^*$ maximizing \eqref{eqn:MLE-column-generating} is added to the set $S$, and the process is repeated. 
Note that \eqref{eqn:MLE-column-generating} has the same form as \eqref{eqn:subgradient-subproblem} and \eqref{eqn:revenue-prediction-robust-constraint}.

We next demonstrate that the MLE problem \eqref{eqn:MLE-frequency} admits a nice interpretation between the empirical estimates $\left\{ p_j \right\}_{j \in [m]}$ and the distributions $\left\{ A_j \lambda \right\}_{j \in [m]}$. To observe this, let us rewrite the objective in \eqref{eqn:MLE-frequency} as
\begin{align*}
\sum_{j \in [m]} \sum_{i \in \cA_j} q_{ij} \log\left(\la a_{ij}, \lambda \ra \right) 
&=  \sum_{j \in [m]} q_j \sum_{i \in \cA_j} p_{ij} \log\left( \la a_{ij}, \lambda \ra  \right)\\
&=  - \sum_{j \in [m]} q_j \underbrace{\sum_{i \in \cA_j} p_{ij} \log\left( \frac{p_{ij}}{\la a_{ij}, \lambda \ra } \right)}_{=  \KL(p_j,A_j \lambda)}  + \underbrace{\sum_{j \in [m]} q_j \sum_{i \in \cA_j} p_{ij} \log(p_{ij})}_{= \text{constant}}
\end{align*}
where $\KL(a,b)$ is the KL divergence between two probability distributions $a$ and $b$.
Hence, \eqref{eqn:MLE-frequency} is equivalent to solving
\begin{equation}\label{eqn:MLE-frequency-alt}
\min_{\lambda} \left\{ \sum_{j \in [m]} q_j \KL(p_j,A_j \lambda) :~  \lambda \in \Delta_{n!} \right\}.
\end{equation}
Thus, by defining $D(A\lambda,p) = \sum_{j \in [m]} q_j \KL(p_j,A_j \lambda)$, we see that the MLE approach is equivalent to \eqref{eqn:revenue-prediction-distance-min} but with a different distance measure $D(\cdot,\cdot)$.

\subsection{Norm-Minimization Approach}\label{sec:norm-min-approach}

As opposed to the approaches outlined in Appendix~\ref{sec:rev-predict}~and~\ref{sec:MLE-approach}, in order to estimate a non-parametric choice model $\lambda$, \citetEC{BertsimasMisic2015} suggest minimizing the $\ell_1$-norm of $p - A \lambda$ by solving 
\begin{equation}\label{eqn:norm-l1-min}
\min_{\lambda} \left\{ \|p - A \lambda\|_1 :~  \lambda \in \Delta_{n!} \right\}.
\end{equation}

In fact, \eqref{eqn:norm-l1-min} can be cast as an LP, but it is still computationally intractable since the dimension of $\lambda$ is $n!$. Similar to \citetEC{vanRyzinVulcano2015}, \citetEC{BertsimasMisic2015} addresses this computational difficulty via a column generation approach. Again, \eqref{eqn:norm-l1-min} is of the same form as \eqref{eqn:revenue-prediction-distance-min} where the distance measure $D(\cdot,\cdot)$ is selected to be $D(A\lambda,p) = \|p - A\lambda\|_1$.
Furthermore, the resulting column generating subproblem is of the form
\begin{equation}\label{eqn:norm-l1-min-column-generating}
\max_{\sigma} \left\{ \sum_{j \in [m]} \sum_{i \in \cA_j} \beta_{ij}(S) a_{ij}(\sigma) - \nu(S) : ~ \sigma \in S_n \right\},
\end{equation}
where $\beta(S)$ and $\nu(S)$ are from the dual solution to solving \eqref{eqn:norm-l1-min} on a subset of columns $\sigma \in S \subset S_n$. Again, this subproblem has the same form as \eqref{eqn:subgradient-subproblem}, \eqref{eqn:revenue-prediction-robust-constraint} and \eqref{eqn:MLE-column-generating}.

%% file: ErrorRates.tex
\section{Convergence Rates for Error Terms $\epsilon^\circ, \epsilon^\bullet$}\label{sec:error-rates}

In Table~\ref{tab:primaldual-JEO-rates}, we state the convergence rate of $\frac{2L}{\Theta_T} \sum_{t \in [T]} \theta_t \|p_t - p\| \to 0$ for different possible rates of $\|p_t - p\| \to 0$, as well as two common choices for $\theta_t$, namely $\theta_t = 1$ and $\theta_t = t$. In Section \ref{sec:JEO-algs}, we discuss the effect of the choice of $\theta$ on the regret bounds for $\widehat{\epsilon}$.
\begin{proposition}\label{prop:primaldual-JEO-rates}
	The convergence rates in Table \ref{tab:primaldual-JEO-rates} hold.
	\begin{table}
		\centering
		\small
		\begin{tabular}{l || c | c}
			rate at which ~$\frac{2L_\Psi}{\Theta_T} \sum_{t \in [T]} \theta_t \|p_t - p\| \to 0$~ &~ $\theta_t = 1$, $\Theta_T = T$~ &~ $\theta_t = t$, $\Theta_T = \frac{T(T+1)}{2}$\\ \hline
			$\|p_t - p\| = O(1/t^r)$, $r \in (0,1)$, & $\sim 1/T^r$ & $\sim 1/T^r$ \\
			$\|p_t - p\| = O(1/t)$ & $\sim \log(T)/T$ & $\sim 1/T$ \\
			$\|p_t - p\| = O(1/t^r)$, $r \in (1,2)$, & $\sim 1/T$ & $\sim 1/T^r$ \\
			$\|p_t - p\| = O(1/t^2)$ & $\sim 1/T$ & $\sim \log(T)/T^2$ \\
			$\|p_t - p\| = O(1/t^r)$, $r > 2$ & $\sim 1/T$ & $\sim 1/T^2$ \\
			$\|p_t - p\| = O(\beta^t)$, $\beta \in (0,1)$ & $\sim 1/T$ & $\sim 1/T^2$
		\end{tabular}
		\caption{Convergence rate of bound for $\epsilon^\circ + \epsilon^\bullet$.}\label{tab:primaldual-JEO-rates}
	\end{table}
\end{proposition}
\proof{Proof of Proposition \ref{prop:primaldual-JEO-rates}}
First, we analyze $S(r,T) := \sum_{t \in [T]} \frac{1}{t^r}$ for $r \neq 1,2$. Observe that
\[ \frac{1}{1-r} \left( \frac{1}{(T+1)^{r-1}} - 1 \right) = \int_{t=1}^{T+1} \frac{1}{t^r} dt \leq S(r,T) \leq \left( 1 + \int_{t=1}^T \frac{1}{t^r} dt \right) = \frac{1}{1-r} \left( \frac{1}{T^{r-1}} - r \right). \]
\begin{itemize}
	\item We consider the case $\theta_t = 1$, $\Theta_T = T$, $\|p_t - p\| = O(1/t^r)$, $r > 0$. In this case we have
	\[ \frac{1}{1-r} \left( \frac{1}{(T+1)^r} - \frac{1}{T+1} \right) \leq \frac{2L}{\Theta_T} \sum_{t \in [T]} \theta_t \|p_t - p\| \sim \frac{1}{T} S(r,T) \leq \frac{1}{1-r} \left( \frac{1}{T^r} - \frac{r}{T} \right). \]
	When $r < 1$, $1-r > 0$ and $1/T = O(1/T^r)$, hence the lower and upper bounds are $\sim 1/T^r$. When $r > 1$, $1-r < 0$ and $1/T^r = O(1/T)$, hence the lower and upper bounds are $\sim 1/T$.
	
	\item Now consider the case $\theta_t = t$, $\Theta_T = 2/(T(T+1))$, $\|p_t - p\| = O(1/t^r)$, $r > 0$. In this case we have
	\begin{align*}
	\frac{2}{2-r} \left( \frac{1}{(T+1)^r} - \frac{1}{(T+1)^2} \right) &\leq \sim \frac{1}{\Theta_T} \sum_{t \in [T]} \theta_t \|p_t - p\|\\
	&\sim \frac{2}{T(T+1)} S(r-1,T) \leq \frac{2}{2-r} \left( \frac{1}{T^{r-1} (T+1)} - \frac{r-1}{T(T+1)} \right).
	\end{align*}
	When $r < 2$, we have $2-r > 0$ and $1/(T(T+1)) = O(1/(T^{r-1}(T+1))) = O(1/T^r)$, hence the lower and upper bounds are $\sim 1/T^r$. When $r > 2$, $2-r < 0$ and $1/(T^{r-1}(T+1)) = O(1/(T(T+1))) = O(1/T^2)$, hence the lower and upper bounds are $\sim 1/T^2$.
	
	\item Now consider the case $\theta_t = 1$, $\Theta_T = T$, $\|p_t - p\| = O(1/t)$. Then
	\[\frac{2L}{\Theta_T} \sum_{t \in [T]} \theta_t \|p_t - p\| \sim \frac{1}{T} \sum_{t \in [T]} \frac{1}{t} \sim \frac{\log(T)}{T}.\]
	
	\item Now consider the case $\theta_t = t$, $\Theta_T = 2/(T(T+1))$, $\|p_t - p\| = O(1/t^2)$. Then
	\[\frac{2L}{\Theta_T} \sum_{t \in [T]} \theta_t \|p_t - p\| \sim \frac{1}{T(T+1)} \sum_{t \in [T]} \frac{1}{t} \sim \frac{\log(T)}{T^2}.\]
\end{itemize}

Finally, consider the case $\|p_t - p\| = O(\beta^t)$ for $\beta \in (0,1)$. When $\theta_t = 1$, $\Theta_T = T$, we have
\[\frac{2L}{\Theta_T} \sum_{t \in [T]} \theta_t \|p_t - p\| = \frac{2L}{\Theta_T} \sum_{t \in [T]} \beta^t = \frac{2L \beta(1-\beta^T)}{T(1-\beta)} \sim 1/T. \]
When $\theta_t = t$, $\Theta_T = 2/(T(T+1))$, we have
\[\frac{4L}{T(T+1)} \sum_{t \in [T]} \theta_t \|p_t - p\| = \frac{4L}{T(T+1)} \sum_{t \in [T]} t \beta^t = \frac{4L \beta(1 - (T+1) \beta^T + T \beta^{T+1})}{T(T+1)(1-\beta)^2} \sim 1/T^2. \]
\Halmos
\endproof

%% file: OCOframework.tex
\section{Online Convex Optimization}\label{sec:oco-intro}
In the standard OCO setting, we are given a convex domain $Z$ and a finite time horizon $T$. In each time period $t \in [T]$, the following takes place:
\begin{itemize}
	\item we make a decision $z_t \in Z$ based on \emph{past} information from time steps $1,\ldots,t-1$ only.
	\item Then, a convex loss function $h_t:Z \to \bbR$ is revealed, we suffer loss $h_t(z_t)$ and get some feedback typically in the form of first-order information $\grad h_t(z_t)$.
\end{itemize}
It is usually assumed that the functions $h_t$ are chosen possibly by an all-powerful adversary that has full knowledge of our learning algorithm---and we know of only the general class of these functions. As such, it is unreasonable to compare the loss of the player across the time horizon to the best possible loss, which would require full knowledge of $h_t$ in advance of choosing $x_t$. Instead, the player's sequence of decisions $z_t$ is evaluated against the best fixed decision in hindsight, and the (average) difference is defined to be the \emph{regret}:
\begin{equation}\label{eqn:regret}
\frac{1}{T} \sum_{t \in [T]} h_t(z_t) - \min_{z \in Z} \frac{1}{T} \sum_{t \in [T]} h_t(x).
\end{equation}
The goal in OCO is to design efficient \emph{regret minimizing} algorithms that generate %
$x_t$ so that the regret tends to zero as $T$ increases. Thus, in OCO we seek algorithms to choose $x_t$ that ensure
\[
\frac{1}{T} \sum_{t \in [T]} h_t(z_t) - \min_{z \in Z} \frac{1}{T} \sum_{t \in [T]} h_t(z) \leq r(T),\quad \lim_{T \to \infty} r(T) = 0,
\]
and the performance of our algorithms is measured by how quickly $r(T)$ tends to $0$.

For our work, we consider two simple modifications to the standard OCO setting: \emph{lookahead decisions} and \emph{weighted regret}. More precisely, lookahead decisions allow for the possibility of choosing $z_t$ \emph{with knowledge of} the loss function $h_t$, while weighted regret modifies \eqref{eqn:regret} to instead be a weighted average
\begin{equation}\label{eqn:weighted-regret2}
\frac{1}{\Theta_T} \sum_{t \in [T]} \theta_t h_t(z_t) - \min_{z \in Z} \frac{1}{\Theta_T} \sum_{t \in [T]} \theta_t h_t(z),
\end{equation}
where (as before) $\{\theta_t\}_{t \in [T]}$ is a collection of non-negative weights, and $\Theta_T = \sum_{t \in [T]} \theta_t$. The reason for examining weighted regret is clear: by appropriately defining $h_t$, the two terms in \eqref{eqn:SP-gap-two-regret} are actually weighted regret terms. We examine lookahead decisions because the $h_t$ that we define to interpret \eqref{eqn:SP-gap-two-regret} as two regret terms, which are $-\Psi(x_t,\cdot;p_t)$ and $\Psi(\cdot,y_t;p_t)$ respectively, actually depend on the decisions $x_t,y_t$ and the data $p_t$, which we have some control over. In particular, in our setting, it is possible to make lookahead decisions for \emph{one} of the regret terms; however, it is not possible to do it for both, because we must choose either $x_t$ before $y_t$ or $x_t$ before $y_t$. The primal oracle algorithms that we introduce in Section \ref{sec:primal-oracle} exactly choose $x_t$ after $y_t$, thus $y_t$ must be chosen in a non-anticipative manner.

Note that \citetEC{HoNguyenKK2016,WangAbernethy2018,AbernethyEtAl2018} have all examined lookahead and weighted regret in OCO before. In fact, \citetEC{HoNguyenKK2016} have used these concepts for \emph{standard} JEO. Our work, however, aims to provide an avenue for using OCO to solve the \emph{saddle point} JEO problem \eqref{eqn:general-SP-JEO}, where one of the domains involved faces an additional high dimensionality challenge. All of these developments are motivated by the problem of dynamic non-parametric estimation of a choice model. In Sections~\ref{sec:learning-choice-model}-\ref{sec:JEO-algs}, we have introduced this problem formally and discuss  the derivation of efficient algorithms for this problem specifically using our general SP JEO framework.

%% file: FOMComparisons.tex
\section{Interpretations of Algorithms from Section \ref{sec:regret-algs}}\label{sec:interpretations}

We discuss interpretations of our methods when applied to $Z = Y$, $z_t = y_t$ for $t \in [T]$ and $h_t(\cdot) = -\Psi(x_t,\cdot;p_t)$. 

\begin{remark}\label{rem:MDequiv}
	For Theorems \ref{thm:OMD} and \ref{thm:s-OMD}, observe that since $x_t = \argmin_{x \in X} \Psi(x,y_t;p_t)$, $\grad h_t(z_t) = -\grad_{y_t} \Psi(x_t,y_t;p_t) = -\grad \left(\min_{x \in X} \Psi(x,y_t;p_t) \right) = -\grad g(y_t;p_t)$. Thus, the gradients we compute are simply the gradients of the (negative of) the dual function $g(y;p_t)$ from \eqref{eqn:general-dual}. In other words, if $p_t = p$ for $t \in [T]$, then using Theorem \ref{thm:OMD} or \ref{thm:s-OMD} to minimize the dual regret $\CR_y$ is simply performing the well-known Mirror Descent algorithm to solve the dual problem \eqref{eqn:general-dual}. When $p_t \neq p$, these update rules %
	are simply performing Mirror Descent on approximate versions of \eqref{eqn:general-dual}, with built-in guarantees on the error when $p_t \neq p$.
	\epr
\end{remark}

\begin{remark}\label{rem:F-Wequiv}
	For Theorem \ref{thm:s-FTL}, let us consider the case when $p_t = p$ for $t \in [T]$. Then $y_t$ are computed as
	\[ y_{t+1} = \argmax_{y \in Y} \left\{ \frac{1}{\Theta_t} \sum_{s \in [t]} \theta_s \Psi(x_s,y; p) \right\}. \]
	If, furthermore, $\Psi(x,y;p)$ is linear in $x$, i.e., it is of the form $\Psi(x,y;p) = \la x, \bm{\Psi}(y;p) \ra - \alpha \omega(y)$, then we can push the sum into the inner product, and since $f(x;p)$ from \eqref{eqn:general-primal} is of the form $f(x;p) = \max_{y \in Y} \Psi(x,y;p) = \max_{y \in Y} \left\{ \la x, \bm{\Psi}(y;p) \ra - \alpha \omega(y) \right\}$, by letting $ \bar{x}_t^\theta = \frac{1}{\Theta_t} \sum_{s \in [t]} \theta_s x_s$ and using the convex envelope theorem we have
	\[ y_{t+1} = \argmax_{y \in Y} \left\{ \left\la \bar{x}_t^\theta, \bm{\Psi}(y; p) \right\ra - \alpha \omega(y) \right\}, \quad \bm{\Psi}(y_{t+1};p) = \grad f(\bar{x}_t^\theta; p). 
	\]
	Now, recalling that
	\begin{align*}
	x_{t+1} &= \argmin_{x \in X} \Psi(x,y_{t+1};p) = \argmin_{x \in X} \left\la x, \bm{\Psi}(y_{t+1};p) \right\ra = \argmin_{x \in X} \la x, \grad f(\bar{x}_t^\theta; p) \ra, \\
	\text{we deduce }~~~ \bar{x}_{t+1}^\theta &= \left( 1 - \gamma_{t+1} \right) \bar{x}_t^\theta + \gamma_{t+1} x_{t+1}, \quad \gamma_{t+1} = \frac{\theta_{t+1}}{\Theta_{t+1}} \in [0,1].
	\end{align*}
	Note that this is exactly a Frank-Wolfe update for the current average point $\bar{x}_t^\theta$ on the primal function $f(x;p)$. Therefore, when $\Psi$ is linear in $x$ and $p_t = p$ for all $t \in [T]$, using Theorem \ref{thm:s-FTL} to minimize $\CR_y$ is equivalent to using the Frank-Wolfe (F-W) algorithm to solve the primal problem \eqref{eqn:general-primal}. See also \citetEC{AbernethyEtAl2018} for an equivalent observation in the case of a particular type of $\Psi$ arising from the convex conjugate of $f$. Thus, within the general context of the JEO problem \eqref{eqn:general-SP-JEO}, we can think of using Theorem \ref{thm:s-FTL} as a generalization of the F-W algorithm to the dynamic setup, with built-in error guarantees for $p_t \neq p$.
	\epr
\end{remark}

%% file: FWcomparison.tex
\section{Rates and Comparison to Frank-Wolfe Methods}\label{sec:FWcomparison}

In our choice model estimation problem, the high-dimensionality challenge of the domain $X$ necessitates the use of projection/prox-free algorithms. In this respect, our developments for  the JEO problem \eqref{eqn:distance-learn-choice-JEO} can be compared against the classical Frank-Wolfe (F-W) algorithm (see \citetEC{Jaggi2013,FreundGrigas2016}), which admits guarantees when using \emph{approximate gradients}. Indeed, under certain assumptions which we will carefully examine, we can think of using $p_t \approx p$ as an approximate gradient method, i.e., $\grad_x D(x,p_t) \approx \grad_x D(x,p)$. Alternatively, \citetEC{DevolderEtAl2014} considers projection-type first-order methods for smooth functions, and provides  guarantees on using approximate gradient oracles within such algorithms. Since projecting onto $X$ defined in \eqref{eqn:X-choice-domain} for the choice model estimation problem is difficult due to the high dimensionality of the domain $X$, we will not discuss the methods of \citetEC{DevolderEtAl2014} and instead focus on the F-W algorithm with approximate gradient oracles in \citetEC{FreundGrigas2016}.

In this appendix, we will do the following.
\begin{itemize}
	\item We first examine the applicability of the standard F-W method due to the issue of the non-smoothness of the objective function for the norm-based distance measures $D(x,p)$.
	\item We then show how online Mirror Descent (MD, Theorem \ref{thm:OMD}) can circumvent the non-smoothness and give the corresponding rates obtainable.
	\item We then present a simple technique to smooth an $\ell_q$-norm for $q \in [2,\infty)$ by squaring it. We examine na\"{i}vely applying the F-W method to solve \eqref{eqn:distance-learn-choice-JEO} and compare the guarantees for F-W with approximate gradients from \citetEC{FreundGrigas2016} to the guarantees from using MD in our framework. We find that the data error terms from the na\"{i}ve F-W method are worse than the ones from using MD in our framework, while the regret bound terms are comparable asymptotically, but involve worse constant factors.
	\item Finally, we present the so-called `Nesterov smoothing' technique for more general norms. We examine na\"{i}vely applying the F-W method to solve \eqref{eqn:distance-learn-choice-JEO} using Nesterov smoothing, and compare the guarantees for F-W with approximate gradients from \citetEC{FreundGrigas2016} to the guarantees from using a F-W algorithm derived from our framework (see Remark \ref{rem:F-Wequiv}). We find that while the regret bound terms are comparable in both settings, the data error terms for the na\"{i}ve F-W method require a certain rate of convergence of $\|p_t - p\| \to 0$ to vanish; this is not a problem for the F-W method derived from our framework.
\end{itemize}

\subsection{Smoothness Requirement for the Frank-Wolfe Algorithm}\label{sec:FW-nonsmooth}
It is known that the F-W method in general \emph{does not} converge on non-smooth objectives; see \citetEC[Example 1]{Nesterov2018} that demonstrates this on a max-type objective function. In the JEO problem \eqref{eqn:distance-learn-choice-JEO} when $D(x,p) = \|x - p\|$, e.g., as in Appendix~\ref{sec:norm-min-approach}, $D(x,p) $ is non-smooth due to the norm. 
In addition, the usual convergence of the F-W algorithm relies on a finite curvature constant assumption that related to the smoothness properties of the function. In particular, the curvature constant $C_D$ of a function $D(x,p)$ of the variable $x$ is defined as
\begin{equation}\label{eqn:curvature-constant}
C_{D} := \sup_{\substack{x,s \in X\\ \alpha \in [0,1]}} \frac{1}{\alpha^2} \left( D((1-\alpha)x+\alpha s,p) - D(x,p) - \alpha \la s-x, \grad_x D(x,p) \ra \right).
\end{equation}
It is well-known that when the function is smooth and the domain is bounded, the associated curvature constant is finite; see e.g., \citetEC[Lemma 7]{Jaggi2013}. Nevertheless, distance measures $D(x,p)$ of interest in the case of non-parametric choice estimation problem, e.g., from Appendix~\ref{sec:learning-choice-model-existing} are non-smooth. 
Moreover, we next show that when $D(x,p)$ is set up based on the norm (see Appendix~\ref{sec:norm-min-approach}) or the KL divergence as in \citetEC{vanRyzinVulcano2015} (see Appendix~\ref{sec:MLE-approach})  the associated curvature constant of $D(x,p)$ is infinite as well.

\begin{proposition}\label{prop:infinite-curvature constant}
Suppose $n>2$. 
For any $q \in [1,\infty]$, the function $D(x,p) = \|x-p\|_q$ has infinite curvature constant \eqref{eqn:curvature-constant} for any $p \in X$. Furthermore, when the MLE based weighted KL divergence is used, i.e., $D(x,p)=\sum_{j \in [m]} w_j \KL(p_j,x_j)$ for any positive weights $w_j$, and $p \in X$ such that $p_{ij} > 0$ for all $i \in A_j$, the curvature constant is infinite.
\end{proposition}
\proof{Proof.} %

We will first show that the curvature constant $C_D$ defined in \eqref{eqn:curvature-constant} of $D(x,p) = \|x-p\|$ is infinite for any $p \in X$. Let us choose $x = p$, reserving the choice of $\alpha \in [0,1]$ and $s \in X$ for later. Then $D(x,p) = 0$, $D((1-\alpha)x + \alpha s,p) = \alpha \|s-p\|$, and the subgradients of $D(x,p)$ are $\left\{ y :~ \|y\|_* \leq 1 \right\}$. Thus, for any selection of subgradient mapping $y(\hat x) \in \grad_x D(\hat x,p)$ we have
\begin{align*}
\frac{1}{\alpha^2} \bigg[ D((1-\alpha)x + \alpha s,p) - D(x,p) - \alpha \la s-x, y(x) \ra \bigg] &= \frac{1}{\alpha^2} \bigg[ \alpha \|s-p\| - \alpha \la s-p, y(x) \ra \bigg]\\
&= \frac{1}{\alpha} \bigg[ \|s-p\| - \la s-p, y(x) \ra \bigg].
\end{align*}
Note that whenever there is a choice $s \in X$ with $\|s-p\| - \la s-p, y(p) \ra > 0$, we can send $\alpha \to 0$ and conclude that the curvature constant $C_D$ is infinite.

To choose the appropriate $s$, we denote the set of subgradients of $\|\cdot\|_q$ at $s-p$ as $G_{\|\cdot\|}(s-p)$. Observe that for a norm $\|\cdot\|$, if $y \in G_{\|\cdot\|}(s-p)$ then $\|y\|_* \leq 1$ and $\la s-p,y\ra = \|s-p\|$. Thus, we need to choose $s \in X$ such that $y(x) \not\in G_{\|\cdot\|}(s-p)$. To do this, we exploit the following property of $\ell_q$-norms. It is simple to check that for $q\in [1,\infty]$ and $y \in G_{\|\cdot\|_q}(s-p)$, we have the property that $y_{ij} > 0 \implies s_{ij} - p_{ij} > 0$. For our selection $y(x)$, first suppose that there exists $i \in \CA_j$ such that $y(x)_{ij} > 0$. Then a ranking $\sigma$ that ranks $i$ last will have $a(\sigma)_{ij} = 0$, so $a(\sigma)_{ij} - p_{ij} \leq 0$ because $p_{ij}\geq 0$. We cannot have $p = a(\sigma)$ for all $(n-1)!$ rankings $\sigma$ that ranks $i$ last (note that $n>2$); hence, there exists one $\sigma$ such that $a(\sigma) \neq p$, and we choose $s = a(\sigma)$. This implies that $y(x)_{ij} > 0$ while $s_{ij}-p_{ij} \leq 0$, hence $y(x) \not\in G_{\|\cdot\|_q}(s-p)$. Now suppose that $y(x)_{ij} \leq 0$ for all item-subset pairs $(i,j)$. If $y(x) = 0$, then the result follows trivially by choosing any $s \neq p$. Suppose now there exists some $y(x) < 0$. It is again simple to check that for $q\in [1,\infty]$ and $y \in G_{\|\cdot\|_q}(s-p)$, we have the property that $y_{ij} < 0 \implies s_{ij} - p_{ij} < 0$. Then a ranking $\sigma$ that ranks $i$ first will have $a(\sigma)_{ij} = 1$, so $a(\sigma)_{ij} - p_{ij} \geq 0$ because $p_{ij}\leq 1$. We cannot have $p = a(\sigma)$ for all $(n-1)!$ rankings $\sigma$ that ranks $i$ first, so there exists one such that $a(\sigma) \neq p$, and we choose $s = a(\sigma)$. This implies that $y(x)_{ij} < 0$ while $s_{ij}-p_{ij} \geq 0$, hence $y(x) \not\in G_{\|\cdot\|_q}(s-p)$. Thus, in all cases for $y(x)$, we can choose the appropriate $s \in X$.

Now consider the weighted KL-divergence $D(x,p) = -\sum_{j \in [m]} w_j \sum_{i \in \CA_j} p_{ij} \log(x_{ij}/p_{ij})$.
We can assume that $p_{ij} > 0$ by simply ignoring terms in the sum for which $p_{ij} = 0$. Choose $x=p$, which ensures that $D(\cdot,p)$ is differentiable at $x$ with $\grad_x D(x,p)_{ij} = -w_j /x_{ij}$. Then we have
\begin{align*}
&\frac{1}{\alpha^2}  \bigg[ D((1-\alpha)x + \alpha s,p) - D(x,p) - \alpha \la s-x, \grad_x D(x,p) \ra  \bigg]\\
&= -\frac{1}{\alpha^2} \sum_{j \in [m]} w_j \sum_{i \in \CA_j} p_{ij} \log\left(1-\alpha + \alpha \frac{s_{ij}}{p_{ij}}\right) + \frac{1}{\alpha} \sum_{j \in [m]} w_j \sum_{i \in \CA_j} \left(\frac{s_{ij}}{p_{ij}} - 1\right).
\end{align*}
Note that the second term is bounded by $\frac{1}{\alpha} \left( \sum_{j \in [m]} w_j \right) (\max_{i,j} 1/p_{ij} - 1)$. Choose $s_{ij} = a(\sigma)$ for any $\sigma \in S_n$. Then there exists some $i,j$ such that $s_{ij} = 0$. Sending $\alpha \to 1$ results in $\log\left(1-\alpha + \alpha \frac{s_{ij}}{p_{ij}}\right) \to \infty$, and the second term is bounded, so the curvature constant $C_D$ is infinite.
\Halmos
\endproof

Note that our framework proposes to handle the non-smoothness of $D(x,p) = \|x - p\|$ due to the norm by defining $\Psi(x,y;p)$ appropriately, and suggests to utilize a regret-minimizing algorithm to bound the dual regret $\CR_y$, from which primal optimality gap bounds can be inferred through \eqref{eqn:SP-gap} (since the dual optimality gap is always non-negative).

\subsection{Basic Convergence Rate Using the Mirror Descent Algorithm}\label{sec:FWcomparison-rates-MD}

We first discuss convergence rates that we can derive within our framework. Since norms are non-smooth, this immediately suggests the utilization of the Mirror Descent algorithm (Theorem \ref{thm:OMD}). Recall that our assumption on $D$ is that it has a representation \eqref{eqn:distance-max-rep}:
\[ D(x,p) = \max_{y \in Y} \left\{ \la B(x-p), y \ra - \alpha \omega(y) \right\} . \]
The first row of Table \ref{tab:summary-distances} shows that any norm $\|\cdot\|$ can be written in the form \eqref{eqn:distance-max-rep} with $\alpha = 0$, $B = I_N$ and $Y = \{y \in \bbR^N : \|y\|_* \leq 1\}$. Define 
\[\tilde{X} := \left\{ x \geq 0 : \sum_{i \in A_j} x_{ij} = 1, \ j \in [m] \right\} \subset \bbR^N, \] 
and note that $X \subset \tilde{X}$, $p_t \in \tilde{X}$ for all $t \in [T]$ (by \eqref{eqn:pq_{ij}}). Thus, assuming that $\omega$ is 1-strongly convex with respect to some possibly different norm $\|\cdot\|_{\omega}$, and defining
\[ \Omega = \max_{\|y\|_* \leq 1} \omega(y) - \min_{\|y\|_* \leq 1} \omega(y), \quad G \geq \max_{x \in X, p' \in \tilde{X}} \|x-p'\|_{\omega,*}, \]
using Theorem \ref{thm:OMD} and Proposition \ref{prop:Y-bounded-continuity} we get the sub-optimality bound
\begin{equation}\label{eqn:MD-rate-bound}
\|\bar{x}^\theta_T - p\| - \min_{x \in X} \|x-p\| \leq \sqrt{\frac{2\Omega G^2}{T}} + \frac{G_Y \|B\|}{T} \sum_{t \in [T]} \|p_t - p\|,
\end{equation}
where $G_Y$ is as defined in Proposition \ref{prop:Y-bounded-continuity}, which in our case will be $1$.

\begin{table}[t!b!]
\centering
\small
\begin{tabular}{ c | c | c | c | c | c | c | c }
$q$ & $Y$ 		  & $B$ & $\omega(y)$ & $\|\cdot\|_{\omega}$ & $\Omega$ & $G$ & $\sqrt{\frac{2\Omega G^2}{T}}$\\ \hline
$1$ & $\{\|y\|_\infty \leq 1\} \subset \bbR^N$ & $I_N$ & $\frac{1}{2} \|y\|_2^2$ & $\|\cdot\|_2$ & $N/2$ & $\sqrt{2m}$ & $\sqrt{\frac{2mN}{T}}$ \\
$1<q<2$ & $\{\|y\|_{\frac{q}{q-1}} \leq 1\} \subset \bbR^N$ & $I_N$ & $\frac{1}{2} \|y\|_2^2$ & $\|\cdot\|_2$ & $N^{\frac{2-q}{q}}/2$ & $\sqrt{2m}$ & $\sqrt{\frac{2mN^{\frac{2-q}{q}}}{T}}$ \\
$2 \leq q < \infty$ & $\{\|y\|_{\frac{q}{q-1}} \leq 1\} \subset \bbR^N$ & $I_N$ & $\frac{q-1}{2} \|y\|_{\frac{q}{q-1}}^2$ & $\|\cdot\|_{\frac{q}{q-1}}$ & $(q-1)/2$ & $(2m)^{1/q}$ & $\sqrt{\frac{(2m)^{2/q}(q-1)}{T}}$ \\
$\infty$ & $\Delta_{2N} \subset \bbR^{2N}$ & $\begin{bmatrix} I_N & -I_N \end{bmatrix}$ & $\sum\limits_{k \in [2N]} y_k \log(y_k)$ & $\|\cdot\|_1$ & $\log(2N)$ & $1$ & $\sqrt{\frac{2\log(2N)}{T}}$
\end{tabular}
\caption{Different $q$-norms and their Mirror Descent constants (Theorem \ref{thm:OMD}).}\label{tab:MD-constants}
\end{table}
Table \ref{tab:MD-constants} gives a summary of possible choices of $\omega$ for different $q$-norms, together with the associated constants $\Omega$, $G$, and the sub-optimality bound. In general, $\omega$ should be chosen so that $\Omega G^2$ is small, and the operation $\argmin_{y \in Y} \left\{ \la z,y \ra + \omega(y) \right\}$ is easy to compute for any $z$.

\subsection{Convergence Guarantees of the Frank-Wolfe Method for Smooth Minimization}\label{sec:FWcomparison-rates-FW}

For a function $D(x,p)$ which is $L$-smooth with respect to a norm $\|\cdot\|$, the standard convergence rate for the F-W method, fixing $p_t = p$, is given in \citetEC[Bound 3.1, Eq. 8]{FreundGrigas2016}:
\[ D(\bar{x}_T,p) - \min_{x \in X} D(x,p) \leq \frac{2L}{T+4} \max_{x,x' \in X} \|x-x'\|. \]
Therefore, if we can use an alternative function $D$ instead of a norm then relate the optimality gap back to the original norm, we can get convergence guarantees for the norm.

In the dynamic setting, however, we need to use an approximate F-W method, since $p_t \neq p$. Essentially, this means that at each iteration $t$, we have an approximate gradient $\grad_x D(x,p_t) \approx \grad_x D(x,p)$, which we require to satisfy
\begin{equation}\label{eqn:approx-grad-FW}
\max_{x,x' \in X} \left| \la \grad_x D(x,p_t) - \grad_x D(x,p), x-x' \ra \right| \leq \delta_t.
\end{equation}
We update according to
\begin{equation}\label{eqn:naive-FW}
x_t = \argmin_{x' \in X} \la \grad_x D(\bar{x}_t,p_t),\; x' \ra, \quad \bar{x}_t = (1-\gamma_t) \bar{x}_{t-1} + \gamma_t x_t.
\end{equation}
When working with approximate gradients, if \eqref{eqn:approx-grad-FW} holds and we choose $\gamma_t$ appropriately, \citetEC[Theorem 5.1, Proposition 5.1]{FreundGrigas2016} provides a convergence rate of
\begin{equation}\label{eqn:approx-naive-FW-rate}
D(\bar{x}_T,p) - \min_{x \in X} D(x,p) \leq \frac{2L}{T+4} \max_{x,x' \in X} \|x-x'\|^2 + \frac{4}{(T+1)(T+2)} \sum_{t \in [T]} (t+1) \delta_t.
\end{equation}

Since norms are non-smooth, we cannot apply the F-W algorithm directly; recall that  the F-W method in general does not converge on non-smooth objectives (see \citetEC[Example 1]{Nesterov2018}). We next discuss two alternative techniques to build smooth approximations to norms.

\subsubsection{Smoothing via Squaring the Norm}\label{app:FW-squared-norms}

A simple way we can smooth a norm is to square it, i.e., $D(x,p) = \|x-p\|^2$. This does not work for arbitrary norms, but it is known that for $\ell_q$-norms, $\frac{1}{2}\|\cdot\|_q^2$ is $(q-1)$-smooth for $q \in [2,\infty)$. 

\paragraph{Convergence Rate with Approximate Gradients.} Using \eqref{eqn:approx-naive-FW-rate}, the na\"{i}ve F-W algorithm using approximate gradients achieves a bound of (with $D_X := \max_{x,x' \in X} \|x-x'\|_q = (2m)^{1/q}$):
\begin{align}
\|\bar{x}_T - p\|_q - \min_{x \in X} \|x-p\|_q &\leq \sqrt{\|\bar{x}_T-p\|_q^2 - \min_{x \in X} \|x-p\|_q^2}\notag\\
&\leq \sqrt{\frac{8 (2m)^{2/q} (q-1)}{T+4}} + \sqrt{\frac{4 (2m)^{1/q}\, (q-1)}{(T+1)(T+2)} \sum_{t \in [T]} (t+1) \|p_t - p\|_q}.\label{eqn:naive-FW-square-smooth-rate}
\end{align}
The first term of \eqref{eqn:naive-FW-square-smooth-rate} is asymptotically comparable to the first term of \eqref{eqn:MD-rate-bound} (third row of Table \ref{tab:MD-constants}), but with worse constants. From Appendix \ref{sec:error-rates}, we see that the second (data error) term of \eqref{eqn:naive-FW-square-smooth-rate} is asymptotically worse than the data error term in \eqref{eqn:MD-rate-bound} for a variety of rates of convergence of $\|p_t - p\| \to 0$, and also has worse constants.

The data error term is derived as follows. The smoothness of $\|x-p\|_q^2$ is derived from the smoothness of $d(z) = \frac{1}{2} \|z\|_q^2$, which is equivalent to
\[ \|\grad d(z) - \grad d(z')\|_* \leq (q-1)\|z-z'\|. \]
Now observe that
\[ \|\grad_x D(x,p_t) - \grad_x D(x,p)\|_{q,*} = \|\grad d(x-p_t) - \grad d(x-p)\|_{q,*} \leq (q-1)\|p_t - p\|_q. \]
Therefore, we have
\begin{align*}
\delta_t &:= \max_{x,x' \in X} \left| \la \grad_x D(x,p_t) - \grad_x D(x,p),\; x-x' \ra \right|\\
&\leq \max_{x,x' \in X} \left\{ \|\grad_x D(x,p_t) - \grad_x D(x,p)\|_{q,*}\; \|x-x'\|_q \right\}\\
&\leq (q-1) \|p_t - p\|_q \max_{x,x' \in X} \|x-x'\|_q = (q-1)  \|p_t - p\|_q D_X.
\end{align*}

The technique of squaring the norm only works for particular $q\in[2,\infty)$. We next propose a more generally applicable method for smoothing.

\subsubsection{Smooth Approximations to General Norms}

We can smooth more general norms by utilizing the `Nesterov smoothing' framework outlined in \citetEC{BeckTeboulle2012}. Interestingly, this also allows us to make use of other algorithms within our framework such as Theorems \ref{thm:s-OMD} and \ref{thm:s-FTL}.

Given a function $g:Z \to \bbR$, we say that $g_{\beta}:Z \to \bbR$ is an $\beta$-approximation of $g$ if there exists $\beta_1,\beta_2 \geq 0$, $\beta_1 + \beta_2 = 1$ such that the following holds:
\begin{equation}\label{eqn:beta-approx-function}
\forall z \in Z, \quad g(z) - \beta_1 \beta \leq g_{\beta}(z) \leq g(z) + \beta_2 \beta.
\end{equation}
\begin{lemma}\label{lemma:smoothing}
Suppose $g,g_\beta$ satisfies \eqref{eqn:beta-approx-function}. Then for any $\bar{z} \in Z$,
\[ g(\bar{z}) - \min_{z \in Z} g(z) \leq g_{L,\beta}(\bar{z}) - \min_{z \in Z} g_{L,\alpha}(z) + \beta. \]
\end{lemma}
\proof{Proof of Lemma \ref{lemma:smoothing}.}
Observe that $g(\bar{z}) \leq g_\beta(\bar{z}) + \beta_1 \alpha$. Second, observe that $\min_{z \in Z} g(z) \geq \min_{z \in Z} g_\beta(z) - \beta_2 \alpha$. Subtracting the appropriate terms gives us the result.
\Halmos
\endproof

Thus, through approximating the non-smooth distance measures $D(x,p)$ with smooth functions, we can use the F-W algorithm and get sub-optimality bounds. For this, we follow the so-called `Nesterov smoothing' technique outlined in \citetEC[Section 4.3]{BeckTeboulle2012}: by setting $\alpha > 0$ instead of $\alpha = 0$ in \eqref{eqn:distance-max-rep} and choosing $\omega$ to be strongly convex, we get a smooth approximation of the norm with the following guarantee.
\begin{lemma}\label{lemma:smooth-distances-approx-quality}
Let
\[ D(x,p) = \max_{y \in Y} \la B(x-p), y \ra, \quad D_{\alpha}(x,p) = \max_{y \in Y} \left\{ \la B(x-p), y \ra - \alpha \omega(y) \right\}. \]
Suppose that $\omega(y) \geq 0$ for all $y \in Y$. Then
\[ D(x,p) - \alpha \max_{y \in Y} \omega(y) \leq D_\alpha(x,p) \leq D(x,p). \]
Similarly, if $\omega(y) \leq 0$ for all $y \in Y$, then
\[ D(x,p) \leq D_\alpha(x,p) \leq D(x,p) + \alpha \max_{y \in Y} |\omega(y)|. \]
\end{lemma}
\proof{Proof of Lemma \ref{lemma:smooth-distances-approx-quality}.}
When $\omega(y) \geq 0$ for all $y \in Y$, we have
\[ \la B(x-p), y \ra \geq \la B(x-p), y \ra - \alpha \omega(y) \geq \la B(x-p),y \ra - \alpha \max_{y' \in Y} \omega(y'). \]
Taking the maximum over $y \in Y$ of all sides gives the first result. The second result is proved similarly.
\Halmos
\endproof

\begin{example}
Let $D(x,p) = \|x-p\| = \max_{y} \left\{\la x-p, y \ra:~:\|y\|_* \leq 1\right\}$. Define the Huber function
\begin{equation}\label{eqn:Huber}
H_\alpha(r) = \begin{cases}
\frac{1}{2 \alpha} r^2 , &r < \alpha\\
|r| - \frac{\alpha}{2}, &r \geq \alpha.
\end{cases}
\end{equation}
Setting $\omega(y) = \frac{1}{2} \|y\|_*^2$, we have
\begin{align*}
D_{\alpha}(x,p) &= \max_{y:\|y\|_* \leq 1} \left\{ \la x-p, y \ra - \frac{\alpha}{2} \|y\|_*^2 \right\} = \max_{\substack{\gamma \in [0,1]\\ y:\|y\|_* = 1}}  \left\{ \gamma \la x-p, y \ra - \frac{\alpha}{2} \gamma^2 \right\} = \max_{\gamma \in [0,1]} \left\{ \gamma \|x-p\| - \frac{\alpha}{2} \gamma^2 \right\}\\
&= \begin{cases}
\frac{1}{2\alpha} \|x-p\|^2, &\|x-p\| \leq \alpha\\
\|x-p\| - \frac{\alpha}{2}, &\|x-p\| > \alpha
\end{cases} \\
&= H_\alpha(\|x-p\|).
\end{align*}
This guarantees that $H_{\alpha}(\|x-p\|)$ is a $\alpha/2$-approximation of $\|x-p\|$. Note that $\omega(y) = \frac{1}{2} \|y\|_*^2$ is not strongly convex in general, which we need in order to use $\omega(y)$ in Theorems \ref{thm:s-OMD} or \ref{thm:s-FTL}. However, if we consider $q$-norms $\|\cdot\| = \|\cdot\|_q$ for $2 \leq q < \infty$, then it is well-known that $\frac{1}{2} \|y\|_*^2 = \frac{1}{2} \|y\|_{\frac{q}{q-1}}^2$ is strongly convex (with respect to $\|\cdot\|_{\frac{q}{q-1}}$) with parameter $(q-1)^{-1}$.

For $q=1$, we can take $\omega(y) = \frac{1}{2} \|y\|_2^2$ to also get an $N/2$-approximation. For $1 < q < 2$, we use the same $\omega$ to get a $N^{\frac{2-q}{q}}/2$-approximation (which is because $\max_{\|y\|_{\frac{q}{q-1}} \leq 1} \|y\|_2 = N^{\frac{1}{2} - \frac{q-1}{q}}$). However, notice that there is no closed form for $D_\alpha$, but this is not a problem for our framework as long as we have the max-type representation \eqref{eqn:distance-max-rep} for $D_\alpha$. For $q=\infty$, we can use the lifted representation of the $\ell_1$-ball into the $2N$-simplex, and set $\omega(y)$ to be the negative entropy, which gets us a $\log(2N)$-approximation. A summary of this is given in Table \ref{tab:smooth-norms}; notice that the choices of $\omega$ and associated constants are quite similar to Table \ref{tab:MD-constants}.
\epr
\end{example}

\begin{table}[t!b!]
\centering
\small
\begin{tabular}{ c | c | c | c | c | c }
$q$ & $Y$ 		  & $B$ & $\omega(y)$ & $\max\limits_{y \in Y} |\omega(y)|$ & $D_{\alpha}(x,p)$ \\ \hline
$1$ & $\{\|y\|_\infty \leq 1\} \subset \bbR^N$ & $I_N$ & $\frac{1}{2} \|y\|_2^2$ & $N/2$ & $\sum\limits_{j \in [m]} \sum\limits_{i \in \CA_j} H_\alpha(x_{ij}-p_{ij})$ \\
$1<q<2$ & $\{\|y\|_{\frac{q}{q-1}} \leq 1\} \subset \bbR^N$ & $I_N$ & $\frac{1}{2} \|y\|_2^2$ & $N^{\frac{2-q}{q}}/2$ & n/a \\
$2 \leq q < \infty$ & $\{\|y\|_{\frac{q}{q-1}} \leq 1\} \subset \bbR^N$ & $I_N$ & $\frac{q-1}{2} \|y\|_{\frac{q}{q-1}}^2$ & $(q-1)/2$ & $H_{\alpha(q-1)}(\|x-p\|_q)$ \\
$\infty$ & $\Delta_{2N} \subset \bbR^{2N}$ & $\begin{bmatrix} I_N & -I_N \end{bmatrix}$ & $\sum\limits_{k \in [2N]} y_k \log(y_k)$ & $\log(2N)$ & $\alpha \log\left( \sum\limits_{j \in [m]} \sum\limits_{i \in \CA_j} 2\cosh\left( \frac{x_{ij}-p_{ij}}{\alpha} \right) \right)$
\end{tabular}
\caption{Different $q$-norms and their smoothings.}\label{tab:smooth-norms}
\end{table}

\paragraph{Convergence Rate with Approximate Gradients.} %
Using a similar argument to the one in Appendix \ref{app:FW-squared-norms}, we can derive the following guarantee on the approximate gradients $\grad_x D_\alpha(x,p_t)$ for functions \eqref{eqn:distance-max-rep}:
\begin{align*}
\delta_t &:= \max_{x,x' \in X} \left| \la \grad_x D_\alpha(x,p_t) - \grad_x D_{\alpha}(x,p), x-x' \ra \right|\\
&\leq \max_{x,x' \in X} \|\grad_x D_\alpha(x,p_t) - \grad_x D_{\alpha}(x,p)\|_* \|x-x'\|\\
&\leq \frac{1}{\alpha} \|p_t - p\| \max_{x,x' \in X} \|x-x'\|.
\end{align*}
Thus, by denoting $D_X := \max_{x,x' \in X} \|x-x'\|$, \eqref{eqn:approx-naive-FW-rate} gives the following guarantee from using the na\"{i}ve F-W method for the JEO problem \eqref{eqn:distance-learn-choice-JEO}:
\[ D_\alpha(\bar{x}_T,p) - \min_{x \in X} D_\alpha(x,p) \leq \frac{2 D_X^2}{\alpha(T+4)} + \frac{4D_X}{\alpha(T+1)(T+2)} \sum_{t \in [T]} (t+1) \|p_t - p\|. \]
Translating this back into a bound on the norm (Lemmas \ref{lemma:smoothing}, \ref{lemma:smooth-distances-approx-quality}), we have
\begin{equation}\label{eqn:naive-FW-nesterov-smooth-rate}
\|\bar{x}_T - p\| - \min_{x \in X} \|x-p\| \leq \frac{2 D_X^2}{\alpha(T+4)} + \alpha \max_{y \in Y} |\omega(y)| + \frac{4D_X}{\alpha(T+1)(T+2)} \sum_{t \in [T]} (t+1) \|p_t - p\|.
\end{equation}
Let us contrast this guarantee with what is achievable in our primal-dual framework by using Theorem \ref{thm:s-FTL} and Proposition \ref{prop:Y-bounded-continuity}. Recall from Remark \ref{rem:F-Wequiv} that this is actually a variant of F-W. By defining 
\[G := \max_{x \in X, p' \in \tilde{X}} \|x-p'\|,\]  
we get the bound
\begin{align}
\|\bar{x}_T^\theta - p\| - \min_{x \in X} \|x-p\| &\leq D_\alpha(\bar{x}_T^\theta,p) - \min_{x \in X} D_\alpha(x,p) + \alpha \max_{y \in Y} |\omega(y)|\notag\\
&\leq \frac{2 G^2}{\alpha(T+1)} + \alpha \max_{y \in Y} |\omega(y)| + \frac{2 G_Y \|B\|}{T(T+1)} \sum_{t \in [T]} t \|p_t - p\|,\label{eqn:framework-FW-nesterov-smooth-rate}
\end{align}
where $G_Y,\|B\|$ are defined as in Proposition \ref{prop:Y-bounded-continuity}; for norms we have $G_Y = 1$, $\|B\| \leq 2$. 
Comparing \eqref{eqn:naive-FW-nesterov-smooth-rate} and \eqref{eqn:framework-FW-nesterov-smooth-rate}, we notice that $G \approx D_X$, hence, the first terms are quite comparable. Also, except for an additive factor $\alpha \max_{y \in Y} |\omega(y)|$, these $O(1/T)$ rates in theory are faster than the $O(1/\sqrt{T})$ rates we get from the non-smooth Mirror Descent. However, if we choose $\alpha$ optimally taking this term into account explicitly, we get rates that are asymptotically comparable with those of the non-smooth Mirror Descent, but with worse constants. We summarize this in Table \ref{tab:smooth-norm-FW-rates}.
\begin{table}[t!b!]
\centering
\small
\begin{tabular}{ c | c | c | c | c | c | c }
$q$ & $Y$ 		  & $\omega(y)$ & $\max\limits_{y \in Y} |\omega(y)|$ & $G=D_X$ & $\alpha^*$ & $\frac{2G^2}{\alpha^*(T+1)} + \alpha^* \max\limits_{y \in Y} |\omega(y)|$ \\ \hline
$1$ & $\{\|y\|_\infty \leq 1\} \subset \bbR^N$ & $\frac{1}{2} \|y\|_2^2$ & $N/2$ & $\sqrt{2m}$ & $\sqrt{\frac{8m}{N(T+1)}}$ & $\sqrt{\frac{8mN}{T+1}}$ \\
$1<q<2$ & $\{\|y\|_{\frac{q}{q-1}} \leq 1\} \subset \bbR^N$ & $\frac{1}{2} \|y\|_2^2$ & $N^{\frac{2-q}{q}}/2$ & $\sqrt{2m}$ & $\sqrt{\frac{8m}{N^{\frac{2-q}{q}}(T+1)}}$ & $\sqrt{\frac{8mN^{\frac{2-q}{q}}}{T+1}}$ \\
$2 \leq q < \infty$ & $\{\|y\|_{\frac{q}{q-1}} \leq 1\} \subset \bbR^N$ & $\frac{q-1}{2} \|y\|_{\frac{q}{q-1}}^2$ & $(q-1)/2$ & $(2m)^{1/q}$ & $\sqrt{\frac{4(2m)^{2/q}}{(q-1)(T+1)}}$ & $\sqrt{\frac{4(2m)^{1/q}(q-1)}{T+1}}$ \\
$\infty$ & $\Delta_{2N} \subset \bbR^{2N}$ & $\sum\limits_{k \in [2N]} y_k \log(y_k)$ & $\log(2N)$ & $1$ & $\sqrt{\frac{2}{\log(2N)(T+1)}}$ & $\sqrt{\frac{8 \log(2N)}{T+1}}$
\end{tabular}
\caption{Different $q$-norms and their (smoothed) convergence rates.}\label{tab:smooth-norm-FW-rates}
\end{table}

The main difference between the bounds for the na\"{i}ve F-W method and the variant from our framework displayed in \eqref{eqn:naive-FW-nesterov-smooth-rate} and \eqref{eqn:framework-FW-nesterov-smooth-rate} respectively is, since our continuity requirement on $\Psi$ (Assumption \ref{ass:JEO-continuity}) is different to \eqref{eqn:approx-grad-FW}, $\alpha$ does not appear in the third data error term of \eqref{eqn:framework-FW-nesterov-smooth-rate}. This is significant because choosing $\alpha \propto 1/\sqrt{T}$ to minimize the first two terms now makes the third data error term $\frac{1}{\sqrt{T}(T+1)} \sum_{t \in [T]} t \|p - p_t\|$, which may diverge even if $\|p - p_t \| \to 0$. Thus, to ensure convergence using the na\"{i}ve F-W method on the smoothed norm with approximate gradients, we need $\|p - p_t\| \to 0$ \emph{sufficiently fast}. In contrast, our primal-dual framework avoids this obstacle completely. We note also that for norms, $G_Y \|B\| \leq 2$, so we can also see from Table \ref{tab:smooth-norm-FW-rates} that the constants in the error terms of \eqref{eqn:framework-FW-nesterov-smooth-rate} are better than \eqref{eqn:naive-FW-nesterov-smooth-rate}.

%% file: Experiments.tex
\section{Computational Study Details and Supplementary Results}\label{sec:computationApp}\label{sec:experiments}

We describe in detail the experimental setup of our computational study. %
All experiments are conducted on a server with 2.8 GHz processor and 64GB memory, using Python 3.6. Gurobi 8.0 (with default Gurobi settings except we limit the number of threads to 2) is used to solve the integer programming subproblems.

\paragraph{Test Instances.} We employ a setup similar to \citetEC[Chapter 4.5.3]{BertsimasMisic2015}. Our ground truth choice model over $n=10$ items (plus one no-choice option) is a mixed MNL model with $K$ segments. Given mixing probabilities $w \in \Delta_K$ and $K$ sets of utilities $\{u_{i,k}\}_{i \in \{0\} \cup [n]}$, $k \in [K]$, the mixed MNL model chooses an item $i \in \CA \subseteq [n]$ with probability
\vspace{-7pt}
\[ \bbP[i \mid \CA] = \sum_{k \in [K]} w_k \frac{u_{i,k}}{u_{0,k} + \sum_{i' \in \CA} u_{i',k}}.\]
\vspace{-13pt}

For each $k \in [K]$, we generate $n+1$ parameters $q_{i,k} \sim U(0,1)$, $i \in \{0\} \cup [n]$ (recall that $0$ denotes the no-choice option present in each subset). The utilities $u_{i,k}$ are then set as follows: four randomly chosen $i \in \{0\} \cup [n]$ are set to $u_{i,k} = Lq_{i,k}$ while the rest are set to $u_{i,k} = q_{i,k}/10$. The mixing probabilities $\{w_k\}_{k \in [K]}$ are chosen randomly from the $(K-1)$-dimensional simplex. We test on 100 randomly generated instances of this ground truth model. 
In the main discussion presented before, we showed results in the dynamic setting for $K=L=5$. In Appendix~\ref{sec:exp-static}, we also provide results in the static setting under various parameter regimes $K\in\{1,5,10\}$ and $L\in\{5,10,100\}$, and thus test the effect of different ground truth models on the conclusions drawn. We observe that in these different ground truth models, the conclusions are in general in line with the ones from $K=L=5$ setting; this supports that the conclusions drawn from the dynamic data experiments with this particular choice of ground truth choice model are likely to be valid for other ground truth choice models as well.

We tested each algorithm on 100 different instances of a ground truth model for each parameter combination.

\paragraph{Distance measures.} In our main discussion, we consider distance measures $D(x,p)$ based on the $\ell_2$-norm. Specifically, we tested $D(x,p)$ from the first and the third rows of Table \ref{tab:summary-distances}, as well as $D(x,p) = \frac{1}{2} \|x-p\|_2^2$, which has max-type representation $\frac{1}{2} \|x-p\|_2^2 = \max_y \left\{ \la x-p, y \ra  - \frac{1}{2} \|y\|_2^2 \right\}$.
In this appendix, we also consider different choices of distance measures $D(\cdot,\cdot)$ arising from $\ell_1$- and $\ell_\infty$-norms. Table \ref{tab:summary-distances-exp} shows the distance measures we used in our experiments for $\ell_1$-, $\ell_2$- and $\ell_\infty$-norms respectively. Note that $\alpha = 0$ recovers the case for the underlying standard norm.
\begin{table}[t!b!]
\centering
\begin{tabular}{ c | c | c | c | c }
$\|\cdot\|$ & $Y$ 		  & $B$ & $\omega(y)$ & $D(x,p)$ \\ \hline
$\|\cdot\|_1$ & $\{\|y\|_\infty \leq 1\} \subset \bbR^N$ & $I_N$ & $\frac{1}{2} \|y\|_2^2$ & $\sum\limits_{j \in [m]} \sum\limits_{i \in A_j} H_\alpha(x_{ij} - p_{ij})$ \\
$\|\cdot\|_2$ & $\{\|y\|_2 \leq 1\} \subset \bbR^N$ & $I_N$ & $\frac{1}{2} \|y\|_2^2$ & $H_\alpha(\|x-p\|_2)$ \\
$\|\cdot\|_\infty$ & $\Delta_{2N} \subset \bbR^{2N}$ & $\begin{bmatrix} I_N & -I_N \end{bmatrix}$ & $\sum\limits_{k \in [2N]} y_k \log(y_k)$ & $\alpha \log\left( \sum\limits_{j \in [m]} \sum\limits_{i \in \CA_j} 2\cosh\left( \frac{x_{ij}-p_{ij}}{\alpha} \right) \right)$
\end{tabular}
\caption{Smoothed norms used in experiments, where $H_\alpha(r) = \frac{1}{2} r^2$ when $r < \alpha$ and $H_\alpha(r) = |r| - \alpha/2$ otherwise is the Huber function.}\label{tab:summary-distances-exp}
\end{table}
Also, when we used a smoothed norm, we tuned the parameter $\alpha$ to minimize the suboptimality gap bound for the \emph{non-smooth} norm. For example, if we have $D_{\alpha}(x,p)$ as the smoothed version of the $\ell_\infty$-norm, then using the F-W algorithm on the static problem, after $T$ iterations we get the bound
\[ D_\alpha(\bar{x}_T^\theta,p) - \min_{x \in X} D_\alpha(x,p) \leq \frac{16}{\alpha T}. \]
However, translating this back into a sub-optimality gap bound for the $\ell_\infty$-norm incurs an additional $\log(2N) \alpha$ additive error, thus we have
\[ \|\bar{x}_T^\theta - p\|_\infty - \min_{x \in X} \|x - p\|_\infty \leq \frac{16}{\alpha T} + \log(2N) \alpha = 8\sqrt{\frac{\log(2N)}{T}}, \]
where the last equality is setting $\alpha = 4/\sqrt{\log(2N)T}$ to minimize the upper bound. See \citetEC{BeckTeboulle2012} for further details on how this is done for other norms.

\paragraph{Algorithm implementation.} We implemented the following solution methods based on our primal-dual framework. For smooth distance measures, i.e., when we have $\alpha > 0$ and $\omega$ is strongly convex in the representation \eqref{eqn:distance-max-rep}, we use the following algorithms:
\begin{itemize}
\item the na\"{i}ve Frank-Wolfe (F-W) algorithm updating
\[ x_{t+1} = \argmin_{x' \in X} \la \grad_x D(\bar{x}_t^\theta,p_t), x' \ra, \quad \bar{x}_{t+1}^\theta = \left( 1- \frac{\theta_{t+1}}{\Theta_{t+1}} \right) \bar{x}_t^\theta + \frac{\theta_{t+1}}{\Theta_{t+1}} x_{t+1}, \]
with $\theta_t = t$. Note that for smooth distances $D$, convergence is guaranteed by our discussion in Appendix \ref{sec:FWcomparison}.
\item The modified F-W algorithm obtained by using Theorem \ref{thm:s-FTL} for the dual updates (see Remark \ref{rem:F-Wequiv}). This is equivalent to updating
\[ x_{t+1} = \argmin_{x \in X} \la \grad_x D(\bar{x}_t^\theta,\bar{p}_t^\theta), x \ra, \quad \bar{p}_t^\theta = \frac{1}{\Theta_t} \sum_{t \in [T]} \theta_t p_t, \quad \bar{x}_{t+1}^\theta = \left( 1- \frac{\theta_{t+1}}{\Theta_{t+1}} \right) \bar{x}_t^\theta + \frac{\theta_{t+1}}{\Theta_{t+1}} x_{t+1}, \]
for $\theta_t = t$. Note that in the static case, when $p_t = p$ for every $t$, this is equivalent to the na\"{i}ve F-W approach.
\item The dual Mirror Descent (MD) algorithm from Theorem \ref{thm:s-OMD}.
\end{itemize}
When $D$ is non-smooth, i.e., it is the norm, we employed the dual MD algorithm from Theorem \ref{thm:OMD}. Note that Theorem \ref{thm:OMD} employed constant step size policies based on constants $\Omega_Y$, $G$ and the maximum iteration count $T$; computing these constants depends on the particular norm and $\omega$ chosen for the domain $Y$ (see Assumption \ref{ass:prox-setup}), but are not difficult to obtain. For each of these methods and norms, we set a maximum iteration limit of $T=10,000$.
 
In the dynamic setting, we use the algorithms outlined above, but only examine the non-parametric estimation model where the distance measure $D(\cdot,\cdot)$ is based on $\ell_2$-norm, and the case of $m=20$ subsets.  
In the static setting, we examine the effect of using different distance measures based on $\ell_1$, $\ell_2$, and $\ell_\infty$ norms as well as the number subsets $m\in\{10,20,50\}$. In Appendix \ref{sec:exp-static}, we provide results for all of these cases for the static setting. The conclusions from these static setup experiments are in line with the base case we discuss here; therefore, we defer the results for distance measures based on other norms and the case of $m \in \{10,50\}$ to Appendix \ref{sec:exp-static}.

Our non-parametric estimation procedure is as follows. We first generate $m$ subsets of $[n]$ of maximum size $\lfloor n/2 \rfloor$ uniformly at random. We append the no-choice option $0$ to all of these (consequently the dimension of the domain $X$ is $\dim(X)=11! \approx 40,000,000$). Using the ground truth model, we compute the $p_{\train}$ vector, where $p_{\train,ij} = \bbP[i \mid A_j]$, and $A_j$ is a subset from our training set. 
In the static setup, we set $p_t = p_{\train}$ at each iteration. In the dynamic setup, we generated a sequence of $p_t \to p_{\train}$, and at each iteration we supply $p_t$ to the algorithm. We initially generate $2000$ random choice observations $(i,j)$, where $A_j$ is one of the training subsets chosen randomly, and $i \in A_j$ is chosen with the probability $p_{\train,ij}$. We then compute $p_1$ using these observations according to \eqref{eqn:pq_{ij}}. For $t \geq 2$, we generate $\kappa \in \bbN$ new observations, then update $p_{t-1}$ with these new observations. We tested various choices of $\kappa$ between $50$ and $1000$.

For both data regimes, we terminate training according to the mean absolute error (MAE), defined as $\MAE(p, p') = \frac{1}{\text{length}(p)} \sum_{i,j} |p_{ij} - p_{ij}'|$, where $\text{length}(p)$ is the length of the vector $p$. In the static setup, we terminate training when $\MAE(\bar{x}_t^\theta, p_{\train}) \leq 0.001$, where $\bar{x}_t^\theta$ is the vector of choice probabilities (for our current estimated model) on the subsets used in training after $t$ iterations. In the dynamic setup, we terminate training when $\MAE(\bar{x}_t^\theta, \bar{p}_t^\theta) \leq 0.001$.

\paragraph{Performance metrics.} We compare the effectiveness of our methods using three criteria: model fit, sparsity, and algorithm efficiency.

To evaluate model fit, we examine the mean absolute error of choice probabilities on subsets generated independently from the training set. Specifically, we generate 100 subsets of $[n]$ of maximum size $\lfloor n/2 \rfloor$ uniformly at random (independently to the training subsets), and append the no-choice option $0$ to each of them. We compute the vector of choice probabilities $p_{\test}$ using our ground truth model. Letting $\bar{x}$ be the choice probabilities on the test subsets computed from the estimated choice model at training termination, we calculate $\MAE(\bar{x}, p_{\test})$.

To evaluate sparsity of our estimated model, we examine the number of different rankings $\sigma$ with positive probability $\lambda(\sigma) > 0$ in our estimated model. Sparsity is very much desired for non-parametric models, since choice probabilities for sparser models can be computed more efficiently.

To evaluate algorithm efficiency, we examine the number of iterations until the termination criterion is reached. While we could have used solution time as another metric for this purpose, we observed in the static setup that solution time is highly correlated  with the number of iterations. Runtimes are affected by how fast the combinatorial subproblem \eqref{eqn:subgradient-subproblem} is solved, but the focus of our work is not on this aspect, hence in our discussions we focused on the number of iterations as a more accurate representation of algorithm efficiency for our purposes.

\subsection{Static Estimation Results}\label{sec:exp-static}

We use the static setup to compare the effect of different parameters (e.g., $K$, $L$, and the norm) on algorithm performance. We compare the F-W algorithm (recall that the na\"{i}ve and modified versions, i.e., $\text{FW}_{\text{na\"ive}}$ and $\text{FW}_{\text{dyn}}$, are equivalent in the static setup, and so we simply refer to it as $\text{F-W}$ in this subsection), as well as the dual MD algorithm for both original non-smoothed norm $\text{MD}_{\text{ns}}$ and its smoothing $\text{MD}_{\text{smth}}$.

\begin{figure}[h!]
\centering
\caption{Performance metrics in the static setup for different $K$, fixing $L=5, m=20$.}\label{fig:compare_L5}
\captionsetup[subfigure]{aboveskip=20pt,belowskip=-5pt}
\includegraphics[page=3,scale=1]{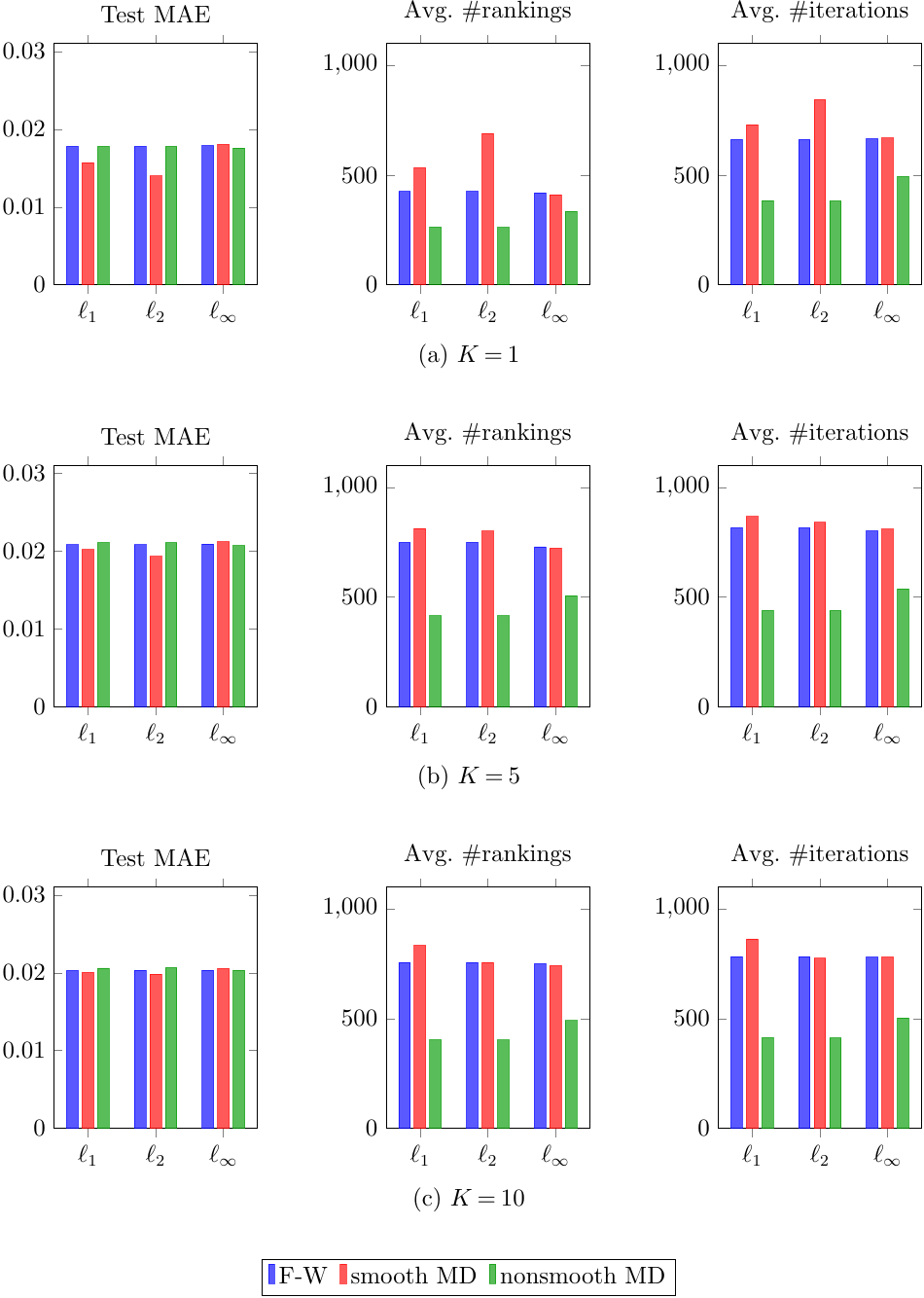}
\end{figure}

\begin{figure}[h!]
\centering
\caption{Performance metrics in the static setup for different $K$, fixing $L=10, m=20$.}\label{fig:compare_L10}
\captionsetup[subfigure]{aboveskip=20pt,belowskip=-5pt}
\includegraphics[page=4,scale=1]{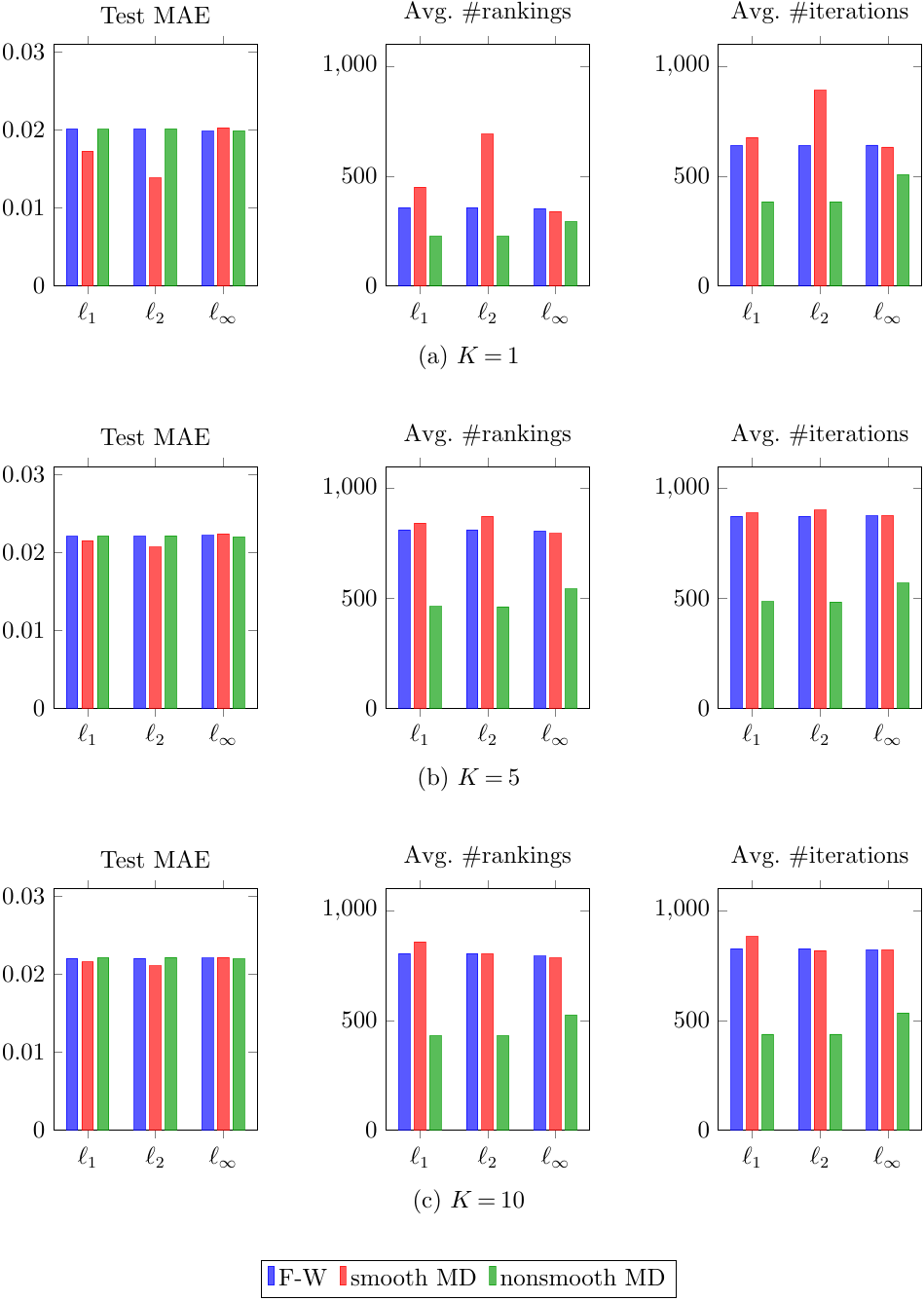}
\end{figure}

\begin{figure}[h!]
\centering
\caption{Performance metrics in the static setup for different $K$, fixing $L=100, m=20$.}\label{fig:compare_L100}
\captionsetup[subfigure]{aboveskip=20pt,belowskip=-5pt}
\includegraphics[page=5,scale=1]{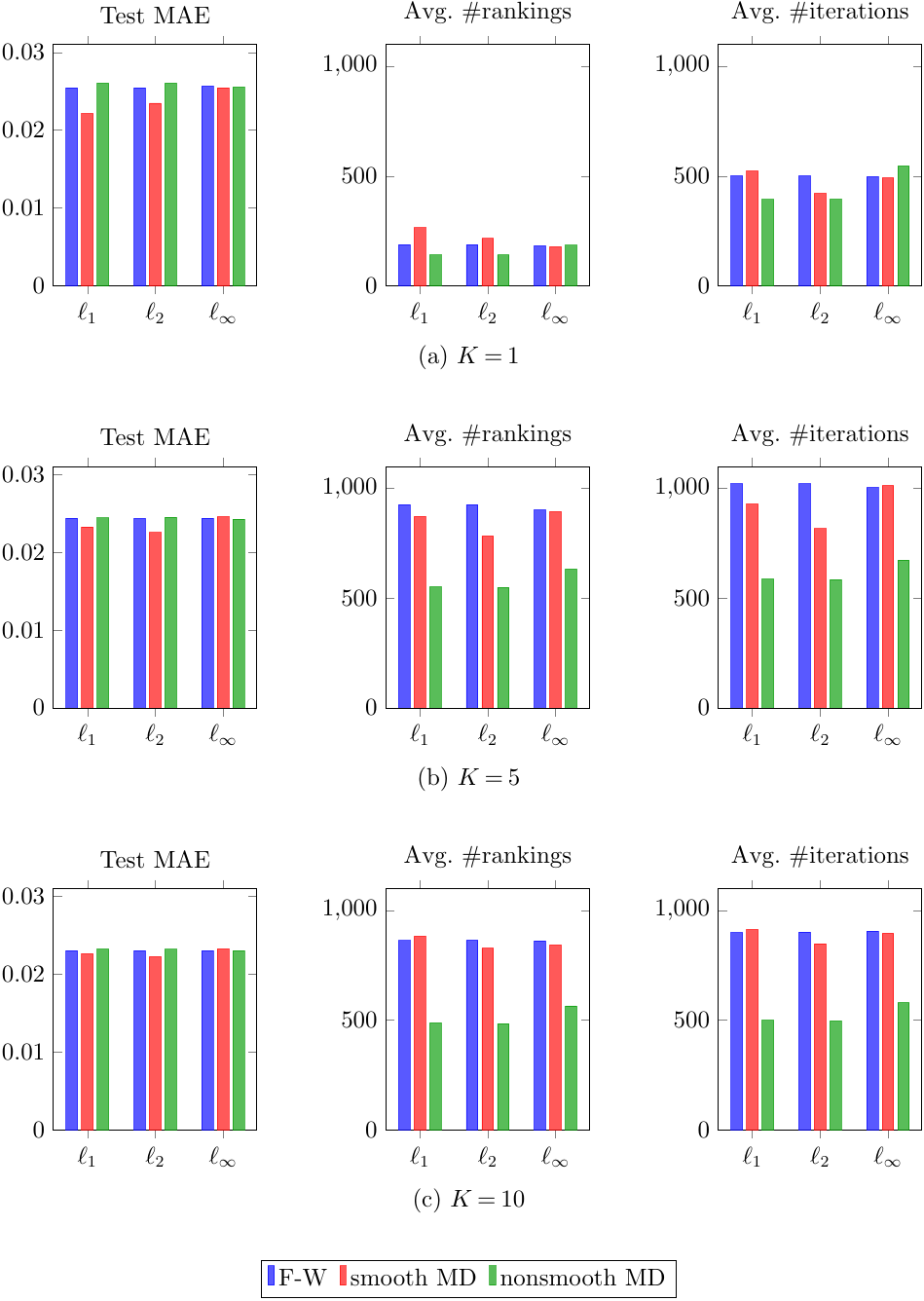}
\end{figure}

\begin{figure}[h!]
\centering
\caption{Performance metrics in the static setup for different $m$, fixing $K=L=5$.}\label{fig:compare_m}
\captionsetup[subfigure]{aboveskip=20pt,belowskip=-5pt}
\includegraphics[page=6,scale=1]{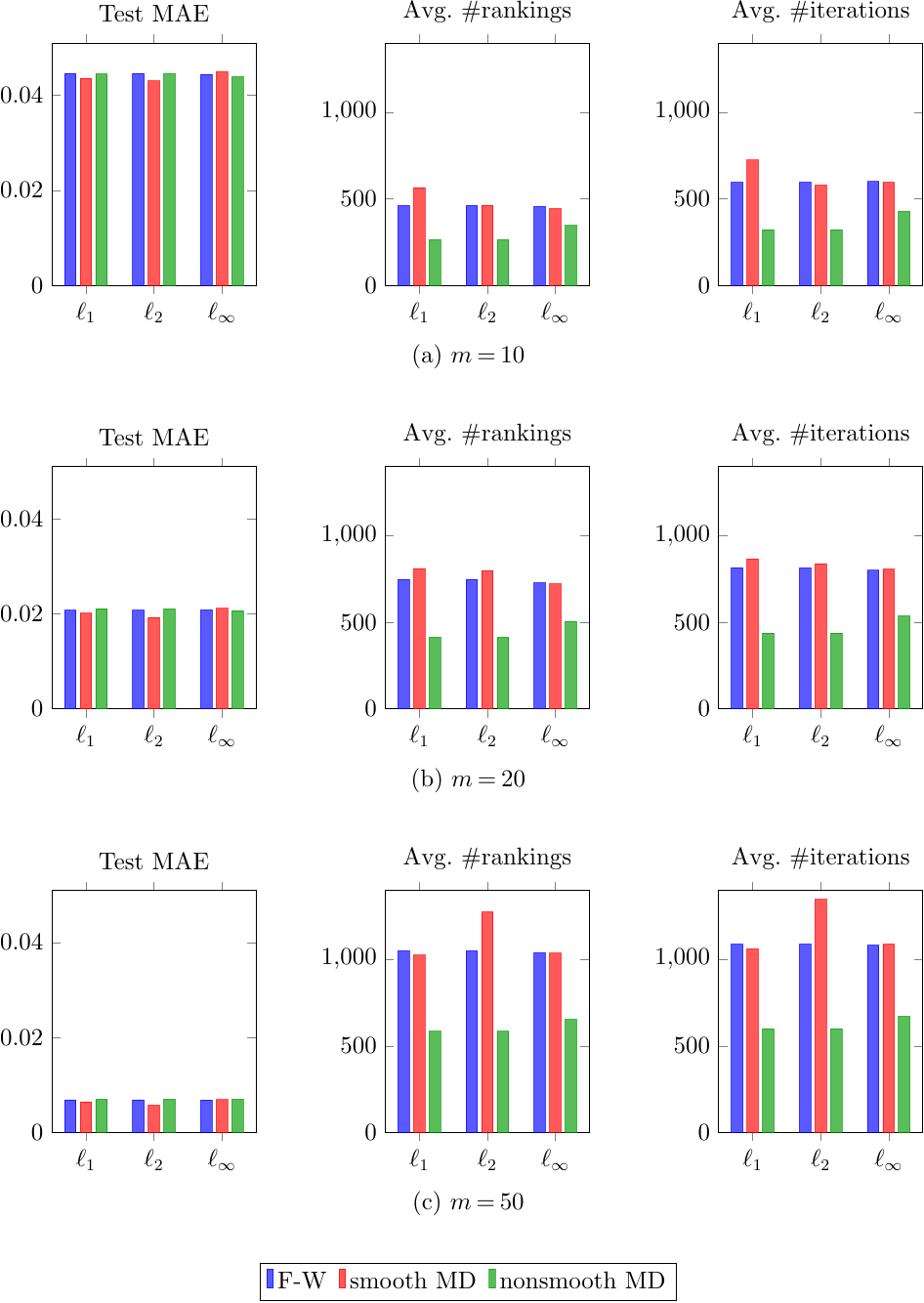}
\end{figure}

\begin{figure}[h!]
\centering
\includegraphics[page=7,scale=1]{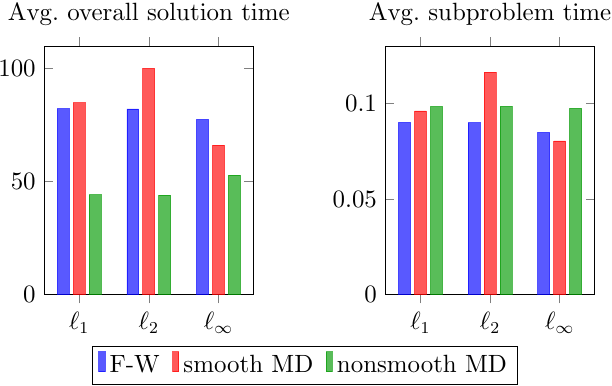}
\caption{Overall solution times and subproblem times per iteration (both in seconds) in the static setup for $m=20$ and $K=L=5$.}
\label{fig:metrics-static-fixedthread}
\end{figure}

Figures \ref{fig:compare_L5}, \ref{fig:compare_L10} and \ref{fig:compare_L100} plot the test MAE, the average number of rankings and the average number of iterations for each of the three solution methods when using different norms, varying $K\in\{1,5,10\}$ and fixing $L$ to respectively $5, 10$ and $100$ while fixing $m=20$. We observe that the three methods have roughly the same test MAE for each of the norms, with perhaps the smoothed dual MD method performing slightly better when $D$ is based on certain norms, but whenever this occurs, it terminates in more iterations and with a denser model. On the other hand, the non-smooth dual MD method clearly learns a sparser model, and clearly terminates in less number of iterations than the other two methods, which are similar in these two metrics. Therefore, we conclude that, regardless of the type of norm used in the estimation procedure, the non-smooth dual MD method is superior in the static setting, since it manages to learn a sparser model more efficiently, while maintaining the same model fit. This conclusion holds for all combinations of $K \in \{1,5,10\}$ and $L = \{5,10,100\}$. In particular, they are consistent with our findings for the dynamic setting results shown in Figure~\ref{fig:metrics-dynamic}.

In Figure \ref{fig:compare_m}, we examine the effect of $m$ by varying $m\in\{10,20,50\}$ while fixing $K=L=5$. We observe that as $m$ increases, the test MAE goes down, but the model sparsity and the number of iterations to convergence increases across all different approaches and norms used. This is as expected, since having more training subsets should allow us to fit better models, but increases the dimension of the choice probability set $X$. Our conclusions regarding the comparison of different approaches remain essentially the same: the non-smooth dual approach still outperforms the others.

Figure \ref{fig:metrics-static-fixedthread} shows the average solution times and the average subproblem times for each method and norm, fixing $m=20$ and $K=L=5$. From Figure \ref{fig:metrics-static-fixedthread} we conclude that the number of iterations and the overall solution time is strongly correlated (contrast it with Figure \ref{fig:compare_m}(b)) and the non-smooth dual approach is still outperforming the other two with respect to overall solution time. We do not believe that the variation in average subproblem time is the result of any inherent property of the methods used or norms. Moreover, the average subproblem solution times are quite small, and thus the variations in subproblem solution times are relatively small.

\subsection{Additional Remarks}

In our numerical experiments, we observed that the average number of rankings and iterations are highly correlated. In fact, the Spearman correlation between these two metrics in the static setting with $K=L=5$ is $\approx 0.922$, thus we conclude that the average number of iterations is a good proxy for model sparsity. This can be seen in the theory: all of our algorithms start with one ranking, and at each iteration they add at most one ranking to the estimated model, which provides an explicit bound on the sparsity of the estimated choice they provide at the end.

Finally, one can argue that the algorithms $\text{FW}_{\text{na\"ive}}$, $\text{FW}_{\text{dyn}}$, and $\text{MD}_{\text{smth}}$ are much simpler to implement than the non-smooth MD algorithm ($\text{MD}_{\text{ns}}$), since there is essentially no parameter tuning aside from computing the smoothing parameters $\alpha$ for $D$ in Table \ref{tab:summary-distances-exp}. On the other hand, $\text{MD}_{\text{ns}}$  algorithm additionally requires tuning the selection of step size, knowing the time horizon $T$, and computing the constants $\Omega_Y$, $G$ (which in turn affects the smoothing parameters). However, for our particular choice model estimation problem, these quantities are quite straightforward to compute, and our analysis and numerical results are based on such `textbook' constant step size policies derived from these, which worked quite well. In terms of performance, we see that the extra sophistication in non-smooth MD ($\text{MD}_{\text{ns}}$) can significantly outperform F-W ($\text{FW}_{\text{na\"ive}}$ or $\text{FW}_{\text{dyn}}$) and smoothed MD ($\text{MD}_{\text{smth}}$).